\documentclass[12pt]{article}
\usepackage{latexsym} \usepackage{epsf}
\usepackage{a4}
\usepackage{amsfonts}
\renewcommand{\theequation}{\thesection.\arabic{equation}}
\newcounter{subequation}[equation]
\makeatletter

\expandafter\let\expandafter\reset@font\csname reset@font\endcsname

\def\subeqnarray{\arraycolsep1pt
    \def\@eqnnum\stepcounter##1{\stepcounter{subequation}%
        {\reset@font\rm(\theequation\alph{subequation})}}
\jot5mm     \eqnarray}

\makeatother

\def\be{\begin{equation}}
\def\ee{\end{equation}}
\def\bea{\begin{eqnarray}}
\def\eea{\end{eqnarray}}
\def\dd{\partial}
\def\half{\frac{1}{2}}

\def\one#1{#1^{\raise5pt\hbox{$\scriptstyle\!\!\!\!1$}}\,{}}
\def\two#1{#1^{\raise5pt\hbox{$\scriptstyle\!\!\!\!2$}}\,{}}

\def\II{\hbox{{1}\kern-.25em\hbox{l}}}
\makeatletter
\def\binrel@#1{\begingroup
  \setboxz@h{\thinmuskip0mu
    \medmuskip\m@ne mu\thickmuskip\@ne mu
    \setbox\tw@\hbox{$#1\m@th$}\kern-\wd\tw@
    ${}#1{}\m@th$}%
  \edef\@tempa{\endgroup\let\noexpand\binrel@@
    \ifdim\wdz@<\z@ \mathbin
    \else\ifdim\wdz@>\z@ \mathrel
    \else \relax\fi\fi}%
  \@tempa
}
\let\binrel@@\relax
\def\overset#1#2{\binrel@{#2}%
  \binrel@@{\mathop{\kern\z@#2}\limits^{#1}}}
\def\underset#1#2{\binrel@{#2}%
  \binrel@@{\mathop{\kern\z@#2}\limits_{#1}}}
\makeatother
\newfont{\bbd}{msbm10 scaled\magstep1}
\def\C{\hbox{\bbd C}}

\def\F{{\mathcal R}}

\def\R{\hbox{\bbd R}}

\def\S{\hbox{\bbd S}}

\def\P{\hbox{\bbd P}}

\def\RR{{\mathcal R}}

\def\RR{{\mathcal R}}

\newtheorem{prop}{Proposition}
\newtheorem{lem}{Lemma}
\begin{document}

\vskip 1cm

\centerline{\LARGE\bf Factorization of the
$\mathrm{R}$-matrix. I.}

\vskip 1cm \centerline{\sc S.E. Derkachov} \vskip 1cm

\centerline{Department of Mathematics, St Petersburg
Technology Institute}\centerline{ St.Petersburg, Russia.}
\centerline{E-mail: {\tt S.Derkachov@pobox.spbu.ru}}

\vskip 2cm

{\bf Abstract.} We study the general rational solution of
the Yang-Baxter equation with the symmetry algebra
$s\ell(3)$. The R-operator acting in the tensor product of
two arbitrary representations of the symmetry algebra can
be represented as the product of the simpler "building
blocks" -- $\RR$-operators. The $\RR$-operators are
constructed explicitly and have simple structure. We
construct in a such way the general rational solution of
the Yang-Baxter equation with the symmetry algebra
$s\ell(3)$. To illustrate the factorization in the simplest
situation we treat also the $s\ell(2)$ case.

\newpage

%%%%%%%%%%%%%%%%%%%%%%%%%%%%%%%%%%%%%%%%%%%%%%%%%%%%%%%%%%%%%%%%%%%%%%%%%%%%%%
{\small \tableofcontents}
\renewcommand{\refname}{References.}
\renewcommand{\thefootnote}{\arabic{footnote}}
\setcounter{footnote}{0}
\setcounter{equation}{0}

\renewcommand{\theequation}{\thesubsection.\arabic{equation}}
\setcounter{equation}{0}

\section{Introduction}
\setcounter{equation}{0}

The Yang-Baxter equation and its solutions play a key role
in the theory of the completely integrable quantum
models~\cite{KS,J,Dr,rep}. The general solution of the
Yang-Baxter equation (R-matrix) is the operator
$\mathbb{R}(u)$ acting in a tensor product
$\mathbf{V}_1\otimes\mathbf{V}_2$ of two linear spaces. The
explicit expression for the R-matrix can be obtained by the
following method~\cite{KRS, MC}. The Yang-Baxter equation
is reduced to the simpler defining equation for the
R-matrix~\cite{KRS}
$$
\R_{12}(u-v)\mathrm{L}_1(u)\mathrm{L}_2(v) =
\mathrm{L}_2(v)\mathrm{L}_1(u)\R_{12}(u-v)
$$
where $\mathrm{L}(u)$ is the Lax operator and by some
conditions the defining equation is reduced to the
recurrence relation for the function of one variable. Let
us consider the $s\ell(2)$-invariant R-matrix~\cite{KRS}
for example. The tensor product of two $s\ell(2)$-modules
has the simple direct sum decomposition
$$
V_{\ell_1}\otimes V_{\ell_2} = \sum_{n=0}^{\infty}
V_{\ell_1+\ell_2+n}
$$
so that one obtains the spectral decomposition of R-matrix
in the following general form
$$
\R_{\ell_1 \ell_2}(u) \sim \sum_{n=0}^{\infty}
\mathrm{R}_n(u)\cdot \P_n
$$
where operator $\P_n$ is the projector on the space
$V_{\ell_1+\ell_2+n}$ in the tensor product
$V_{\ell_1}\otimes V_{\ell_2}$. The defining equation
results in the recurrence relation for the function
$\mathrm{R}_n(u)$
$$
(-u+\ell_1+\ell_2+n)\cdot\mathrm{R}_{n+1}(u)=
-(u+\ell_1+\ell_2+n)\cdot\mathrm{R}_{n}(u)
$$
which has the solution~\cite{rep}
\be \mathrm{R}_n(u) =
(-1)^{n} \frac{\Gamma\left(u+\ell_1+\ell_2+n\right)}
{\Gamma\left(-u+\ell_1+\ell_2+n\right)}. \label{spec} \ee
This method is generally applicable provided the direct sum
decomposition of the tensor product
$\mathbf{V}_1\otimes\mathbf{V}_2$ has no
multiplicities~\cite{MC}. In generic situation for the
algebras of higher rank ($s\ell(3)$ is the first nontrivial
example) one obtains the system of complicated recurrence
relations for the function of several
variables~\cite{KRS,AK,KR} which is difficult to solve.

We suggest the natural factorized expression for the
general R-matrix. It can be represented as the product of
the simple "building blocks" -- $\RR$-operators. The main
idea is very simple and can be illustrated on the
$s\ell(2)$-example. The Lax operator depends on two
parameters: spin of representation $\ell$ and spectral
parameter $u$. It is useful to switch to another two
parameters $u_1 = u+\ell$ and $u_2 = u-\ell$ and extract
the operator of permutation $\P_{12}$ from the R-matrix
$\R_{12}=\P_{12}\check{\R}_{12}$. The defining equation for
the operator $\check{\R}_{12}$ has the form
$$
\check{\R}_{12}\cdot
\mathrm{L}_1(u_1,u_2)\mathrm{L}_2(v_1,v_2) =
\mathrm{L}_1(v_1,v_2)\mathrm{L}_2(u_1,u_2)\cdot
\check{\R}_{12}.
$$
The operator $\check{\R}_{12}$ interchanges $u_1$ with
$v_1$ and $u_2$ with $v_2$ in the product of two
Lax-operators. Let us perform this operation in two stages.
In the first stage we interchange the parameters $u_1$ with
$v_1$ only. The parameters $u_2$ and $v_2$ remain the same.
In this way one obtains the natural defining equation for
the $\RR_1$-operator
$$
\RR_1\cdot \mathrm{L}_1(u_1,u_2)\mathrm{L}_2(v_1,v_2) =
\mathrm{L}_1(v_1,u_2)\mathrm{L}_2(u_1,v_2)\cdot\RR_{1}
$$
In the second stage we interchange $u_2$ with $v_2$ but the
parameters $u_1$ and $v_1$ remain the same. The defining
equation for the $\RR_2$-operator is
$$
\RR_2\cdot \mathrm{L}_1(u_1,u_2)\mathrm{L}_2(v_1,v_2) =
\mathrm{L}_1(u_1,v_2)\mathrm{L}_2(v_1,u_2)\cdot\RR_{2}
$$
These equations appear much simpler then initial defining
equation for the R-operator and the solution can be
obtained in a closed form. Finally we construct the
"composite object" R-matrix from the simplest "building
blocks" -- $\RR$-operators $\R_{12} = \P_{12}\RR_1 \RR_2$.
Note this factorization is different from the ones used by
V.Drinfeld~\cite{Dr1} and the $\RR$-operator is different
from the Drinfeld twist $\mathcal{F}$~\cite{Dr1,Twist}. In
the present paper we shall consider two examples of such
factorization. As an illustrative example we derive the
factorized expression for the R-matrix with symmetry
algebra $s\ell(2)$. Next we work out in details the first
nontrivial example of the R-matrix with symmetry algebra
$s\ell(3)$. It seems that the whole construction can be
generalized to the case of the R-matrix with symmetry
algebra $s\ell(n)$. There exist operators $\RR_1 , \RR_2 ,
\cdots \RR_n$ and the general R-matrix can be expressed in
the factorized form $\R_{12} = \P_{12}\RR_1 \RR_2 \cdots
\RR_n$.

The presentation is organized as follows. In Section 2 we
collect the standard facts about the algebra~$s\ell(2)$ and
its representations. We represent the lowest weight modules
by polynomials in one variable~($z$) and the
$s\ell(2)$-generators as first order differential
operators. We derive the defining relation for the general
R-matrix, i.e. the solution of the Yang-Baxter equation
acting on tensor products of two arbitrary representations,
the elements of which are polynomials in~$z_1$ and~$z_2$.
Next we introduce the natural defining equations for the
$\F$-operators and show that the general R-matrix can be
represented as the product of such much more simple
operators. In the Section 3 we follow the same strategy for
the algebra~$s\ell(3)$. We represent the~$s\ell(3)$ lowest
weight modules by polynomials in three variables~($x , y,
 z$) and the $s\ell(3)$-generators as first order
differential operators. Next we derive the defining
relation for the general R-matrix and show how it can be
solved using the $\F$-operators. Finally, in Section 4 we
summarize.

\section{The general $sl(2)$-invariant R-matrix}
\setcounter{equation}{0}

In this section we consider the simplest situation when the
symmetry algebra of the Yang-Baxter equation is $s\ell(2)$.
We serve it as the illustrative example and present all
calculations in the great details.

\subsection{$sl(2)$ lowest weight modules}
\setcounter{equation}{0}

In this, preparatory, section we collect some facts about
$s\ell(2)$ lowest weight modules and fix the notations. The
Lie algebra $s\ell(2)$ has three generators $\mathbf{S}\ ,\
\mathbf{S}_{\pm}$
$$
\left[\mathbf{S},\mathbf{S}_{\pm}\right] = \pm
\mathbf{S}_{\pm} \ ,\
\left[\mathbf{S}_{+},\mathbf{S}_{-}\right] = 2 \mathbf{S}\
$$
the central element (Casimir operator $\mathbf{C}_2$) being
$$
\mathbf{C}_2 = \mathbf{S}^2-\mathbf{S} +
\mathbf{S}_{+}\mathbf{S}_{-} \ ,\
\left[\mathbf{C}_2,\mathbf{S}_{\pm}\right] =
\left[\mathbf{C}_2,\mathbf{S}\right] = 0.
$$
The Verma module $\mathbf{V}_{\ell}$ is the generic lowest
weight $s\ell(2)$-module with the lowest weight $\ell \in
\C$ and Casimir $\mathrm{C}_2 = \ell(\ell-1)$. As a linear
space $\mathbf{V}_{\ell}$ is spanned by the basis
$\left\{\mathbf{v}_k\right\}_{k=0}^{\infty}$
$$
\mathbf{v}_k = \mathbf{S}_{+}^k\mathbf{v}_{0} \ ,\
\mathbf{S}\mathbf{v}_k = (\ell+k)\mathbf{v}_k \ ,\
\mathbf{S}_{-}\mathbf{v}_k = -k(2\ell+k-1)\mathbf{v}_k
$$
where the vector $\mathbf{v}_0$ is the lowest weight
vector: $\mathbf{S}_{-}\mathbf{v}_0 = 0\ ,\
\mathbf{S}\mathbf{v}_0 = \ell \mathbf{v}_0$. The module
$\mathbf{V}_{\ell}$ is irreducible, except for $\ell = -
\frac{n}{2}$ where $ n \in \{0, 1, 2, 3 \cdots\}$, when
there exists an $(n+1)$-dimensional invariant subspace
$\mathbf{V}_{n}\subset\mathbf{V}_{\ell}$ spanned by
$\left\{\mathbf{v}_k\right\}_{k=0}^{n}$. We shall
extensively use the representation $\mathrm{V}_{\ell}$ of
$s\ell(2)$ in the infinite-dimensional space $\C[z]$ of
polynomials in variable $z$ with the standard monomial
basis $\left\{z^k\right\}_{k=0}^{\infty}$ and lowest weight
vector $v_0 = 1$. The action of $s\ell(2)$ in
$\mathrm{V}_{\ell}$ is given by the first-order
differential operators: \be \mathrm{S} = z\dd + \ell \ ,\
\mathrm{S}_{-} = -\dd\ ,\ \mathrm{S}_{+} = z^2\dd + 2\ell z
\label{diff} \ee or the group-like elements(global
transformations) \be \lambda^{\mathrm{S}}\Phi(z) =
\lambda^{\ell}\Phi(\lambda z)\ ,\ \mathrm{e}^{\lambda
\mathrm{S}_{-}}\Phi(z) =\Phi\left(z-\lambda\right)\ ,\
\mathrm{e}^{\lambda \mathrm{S}_{+}}\Phi(z) = (1-\lambda
z)^{-2\ell}\Phi\left(\frac{z}{1-\lambda z}\right)
\label{glob} \ee The generating function for the basis
vectors of Verma module $\mathbf{V}_{\ell}$
$$
\mathrm{e}^{\lambda \mathbf{S}_{+}}\cdot \mathbf{v}_0 =
\sum_{k=0}^{\infty} \frac{\lambda^k}{k!}\cdot
\mathbf{S}_{+}^{k} \mathbf{v}_{0} = \sum_{k=0}^{\infty}
\frac{\lambda^k}{k!}\cdot \mathbf{v}_{k}
$$
can be calculated in closed form in the functional
representation $\mathrm{V}_{\ell}$ using~(\ref{glob})
$$
\mathrm{e}^{\lambda \mathrm{S}_{+}}\cdot 1 = (1-\lambda
z)^{-2\ell} = \sum_{k=0}^{\infty} \frac{\lambda^k}{k!}\cdot
\left(2\ell\right)_k z^{k} \ ;\ \left(2\ell\right)_k \equiv
\frac{\Gamma(2\ell+k)}{\Gamma(2\ell)}
$$
This expression clearly shows that for generic $\ell\neq
-\frac{n}{2}$ the module $\mathrm{V}_{\ell}$ is an
irreducible $s\ell(2)$-module isomorphic to
$\mathbf{V}_{\ell}$, the isomorphism being given by
$\mathbf{v}_k \leftrightarrow \left(2\ell\right)_k\cdot
z^k$. For $\ell = -\frac{n}{2}$ where $ n \in \{0, 1, 2, 3
\cdots\}$, we have the finite sum instead of infinite
series
$$
\mathrm{e}^{\lambda \mathrm{S}_{+}}\cdot 1 = (1-\lambda
z)^{n} = \sum_{k=0}^{n} (-1)^k C^n_k \lambda^k z^k\ ;\
C^{n}_k = \frac{n!}{k!(n-k)!}
$$
so that there exists an invariant subspace
$\mathrm{V}_{n}\subset\mathrm{V}_{\ell}$ spanned by
$\left\{z^k\right\}_{k=0}^{n}$ which is isomorphic to
$\mathbf{V}_{n}$. For $\ell = -\frac{1}{2}$ one obtains the
two-dimensional invariant subspace $\mathrm{V}_{1}\sim
\C^2$ and the matrices of operators $\mathrm{S}\ ,\
\mathrm{S}^{\pm}$ in the basis
$$
\mathbf{e}_1 = \mathrm{S}_{+}\cdot 1 = -z ,\ \mathbf{e}_2 =
1
$$ have the standard form of generators
$\mathbf{s} , \mathbf{s}_{\pm}$ in the fundamental
representation
\begin{equation}
\mathbf{s}_{+} = \left(\begin{array}{cc}
0 & 1 \\
0 & 0\end{array}\right)\ \ ,\ \mathbf{s}_{-} =
\left(\begin{array}{cc}
0 & 0 \\
1 & 0\end{array}\right)\ \ ,\ \mathbf{s} =
\half\cdot\left(\begin{array}{cc}
1 & 0 \\
0 & -1\end{array}\right)\ \label{fun}
\end{equation}

\subsection{Yang-Baxter equation and Lax operator}
\setcounter{equation}{0}

Let $\mathrm{V}_{\ell_1}$, $\mathrm{V}_{\ell_2}$ and
$\mathrm{V}_{\ell_3}$ be lowest weight $s\ell(2)$-modules
and consider the three operators $\R_{\ell_i\ell_j}(u)$
which are acting in $V_{\ell_i}\otimes V_{\ell_j}$. The
Yang-Baxter equation is the following three term relation
\be \R_{\ell_1\ell_2}(u-v)\R_{\ell_1\ell_3}(u)
\R_{\ell_2\ell_3}(v)=
\R_{\ell_2\ell_3}(v)\R_{\ell_1\ell_3}(u)
\R_{\ell_1\ell_2}(u-v) \label{YB} \ee We seek the general
$s\ell(2)$-invariant solution $\R_{\ell_1\ell_2}(u)$ of
this equation. The natural way is to start from the
simplest solutions of Yang-Baxter equation and derive the
defining equation for the general $\R$-operator~\cite{KRS}.
First we put $\ell_1 = \ell_2 = \ell_3 = -\frac{1}{2}$
in~(\ref{YB}), consider the restriction to the invariant
subspace $\C^2\otimes\C^2\otimes\C^2$ and obtain the
equation
$$ \R_{1 2}(u-v)\R_{1 3}(u) \R_{2 3}(v)= \R_{2
3}(v)\R_{1 3}(u) \R_{1 2}(u-v) $$ where the operator
$\R_{12}(u)$ acts on the first and second copy of $\C^2$ in
the tensor product $\C^2\otimes\C^2\otimes\C^2$ and
similarly for the other $\R$-operators. The solution is the
Yang's $s\ell(2)$-invariant $\mathrm{R}$-matrix~\cite{Y}
$$
\R_{12}(u) = u + \mathrm{P}_{12}
$$
where $\mathrm{P}_{12}$ is the permutation operator in
$\C^2\otimes\C^2$. Next we choose $\ell_1 = \ell\ ;\ \ell_2
= \ell_3 = -\half$ and consider the restriction to the
invariant subspace
$\mathrm{V}_{\ell}\otimes\C^2\otimes\C^2$. The restriction
of the operator $\R_{\ell , -\half}(u)$ to the space
$\mathrm{V}_{\ell}\otimes\C^2$ coincides up to
normalization and the shift of the spectral parameter with
the fundamental Lax-operator
$$
\mathrm{L}(u): \mathrm{V}_{\ell}\otimes\C^2 \to
\mathrm{V}_{\ell}\otimes\C^2
$$
and the Yang-Baxter equation defines the commutation
relations for the Lax-operators~\cite{rep}
$$
\mathrm{L}^{(1)}(u-v)\mathrm{L}^{(2)}(u) \R_{1 2}(v)= \R_{1
2}(v)\mathrm{L}^{(2)}(u) \mathrm{L}^{(1)}(u-v)
$$ where $\mathrm{L}^{(1)}(u)$ is the operator which acts
nontrivially on the space $\mathrm{V}_{\ell}$ and the first
copy of $\C^2$ in the tensor product
$\mathrm{V}_{\ell}\otimes\C^2\otimes\C^2$ and
$\mathrm{L}^{(2)}(u)$ is the operator which acts
nontrivially on the $\mathrm{V}_{\ell}$ and the second copy
of $\C^2$. The solution is given (up to additive constant)
by the Casimir operator $\mathbf{C}_2$ for the tensor
product of representations
$\mathrm{V}_{\ell}\otimes\C^2$\cite{KRS}
$$
\mathrm{L}(u) \equiv u  + 2\mathrm{S}\otimes
\mathbf{s}+\mathrm{S}_{-}\otimes
\mathbf{s}_{+}+\mathrm{S}_{+}\otimes\mathbf{s}_{-} =
\left(\begin{array}{cc}
u+\mathrm{S} & \mathrm{S}_{-} \\
\mathrm{S}_{+} & u-\mathrm{S}
   \end{array}\right) = \left(\begin{array}{cc}
u+\ell+z\dd & -\dd \\
z^2\dd + 2\ell z& u-\ell-z\dd
   \end{array}\right)
$$
where $\mathbf{s}, \mathbf{s}^{\pm}$ are the generators in
the fundamental representation~(\ref{fun}) and $\mathrm{S},
\mathrm{S}_{\pm}$ are generators~(\ref{diff}) in the
generic representation~$\mathrm{V}_{\ell}$. The Lax
operator acts in the space $\C[z]\otimes\C^2$ and despite
the compact notation $\mathrm{L}(u)$ depends really on two
parameters: spin $\ell$ and spectral parameter $u$. We
shall use extensively the parametrization $u_{1}\equiv
u+\ell, u_{2}\equiv u-\ell$ and show all parameters
explicitly. There exists very useful factorized
representation for the $\mathrm{L}$-operator~\cite{Skl1}
\begin{equation}
\mathrm{L}(u_{1},u_{2}) \equiv \left(\begin{array}{cc}
u_{1} +z\dd & -\dd \\
z^2\dd + (u_{1} -u_{2}) z& u_{2} -z\dd
   \end{array}\right) = \left(\begin{array}{cc}
1 & 0 \\
z& 1\end{array}\right)\ \left(\begin{array}{cc}
u_{1} -1 & -\dd \\
0& u_{2}\end{array}\right)\ \left(\begin{array}{cc}
1 & 0 \\
-z& 1\end{array}\right) \label{factor}
\end{equation}
The $\mathrm{L}$-operator is $s\ell(2)$ invariant by
construction and as consequence we obtain the useful
equality \be \left(\begin{array}{cc}
1 & 0 \\
-\lambda & 1\end{array}\right)\cdot\mathrm{L}(u)\cdot
\left(\begin{array}{cc}
1 & 0 \\
\lambda & 1\end{array}\right)
=\mathrm{e}^{\lambda\mathrm{S}_{-}}\cdot\mathrm{L}(u)\cdot
\mathrm{e}^{-\lambda\mathrm{S}_{-}}=
\mathrm{e}^{-\lambda\partial}\cdot\mathrm{L}(u)\cdot
\mathrm{e}^{\lambda\partial} \label{sl2} \ee Indeed, we
have
$$
\left(\mathbf{s}_{-}+\mathrm{S}_{-}\right)\cdot\mathrm{L}(u)
=
\mathrm{L}(u)\cdot\left(\mathbf{s}_{-}+\mathrm{S}_{-}\right)
\Rightarrow
\mathrm{e}^{-\lambda\mathbf{s}_{-}}\cdot\mathrm{L}(u)\cdot
\mathrm{e}^{\lambda\mathbf{s}_{-}} =
\mathrm{e}^{\lambda\mathrm{S}_{-}}\cdot\mathrm{L}(u)\cdot
\mathrm{e}^{-\lambda\mathrm{S}_{-}}\ ;\
\mathrm{e}^{\lambda\mathbf{s}_{-}} =
\left(\begin{array}{cc}
1 & 0 \\
\lambda & 1\end{array}\right)
$$
Finally we put $\ell_3 = -\frac{1}{2}$ in~(\ref{YB}) and
consider the restriction on the invariant subspace
$\mathrm{V}_{\ell_1}\otimes\mathrm{V}_{\ell_2}\otimes\C^2$.
In this way one obtains the defining equation for the
operator $\R_{\ell_1\ell_2}(u)$~\cite{KRS,Skl}
$$
\R_{\ell_1\ell_2}(u-v)\mathrm{L}_1(u) \mathrm{L}_2(v) =
\mathrm{L}_2(v)\mathrm{L}_1(u) \R_{\ell_1\ell_2}(u-v)
$$
The operator $\mathrm{L}_k$ acts nontrivially on the tensor
product $\mathrm{V}_{\ell_k}\otimes\C^2$ which is
isomorphic to $\C[z_k]\otimes\C^2$ and the operator
$\R_{\ell_1\ell_2}(u)$ acts nontrivially on the tensor
product $\mathrm{V}_{\ell_1}\otimes\mathrm{V}_{\ell_2}$
which is isomorphic to $\C[z_1]\otimes\C[z_2] =
\C[z_1,z_2]$.

\subsection{The general R-matrix}
\setcounter{equation}{0}

Now we are going to the solution of the defining equation
for the general R-matrix. It is useful to extract the
operator of permutation $\P_{12}$
$$
\P_{12}\Psi(z_1,z_2) = \Psi(z_2,z_1)\ ;\
\Psi(z_1,z_2)\in\C[z_1,z_2]
$$
from the $\R$-operator $\R_{\ell_1\ell_2}(u) =
\P_{12}\check{\R}_{\ell_1\ell_2}(u)$ and solve the defining
equation for the $\check{\R}$-operator
$$
\check{\R}(u_{1},u_2|v_{1},v_2)
\mathrm{L}_{1}(u_1,u_2)\mathrm{L}_{2}(v_1,v_2)=
\mathrm{L}_{1}(v_1,v_2)\mathrm{L}_{2}(u_1,u_2)
\check{\R}(u_{1},u_2|v_{1},v_2).
$$
To avoid misunderstanding we present this equation in
explicit form
$$
\check{\R}(u_{1},u_2|v_{1},v_2) \left(\begin{array}{cc}
u_{1} +z_1\dd_1 & -\dd_1 \\
z_1^2\dd_1 + (u_{1} -u_{2}) z_1 & u_{2} -z_1\dd_1
   \end{array}\right)\cdot\left(\begin{array}{cc}
v_{1} +z_2\dd_2 & -\dd_2 \\
z_2^2\dd_2 + (v_{1} -v_{2}) z_2& v_{2} -z_2\dd_2
   \end{array}\right)\ = $$$$
 =  \left(\begin{array}{cc}
v_{1} +z_1\dd_1 & -\dd_1 \\
z_1^2\dd_1 + (v_{1} -v_{2}) z_1 & v_{2} -z_1\dd_1
   \end{array}\right)\cdot\left(\begin{array}{cc}
u_{1} +z_2\dd_2 & -\dd_2 \\
z_2^2\dd_2 + (u_{1} -u_{2}) z_2& u_{2} -z_2\dd_2
   \end{array}\right)\
   \check{\R}(u_{1},u_2|v_{1},v_2)
$$
where $u_1 = u+\ell_1\ ,\ u_2 = u-\ell_1\ ,\ v_1 =
v+\ell_2\ ,\ v_2 = v - \ell_1$. The $\check{\R}$-operator
can be factorized on the product of the simpler "elementary
building blocks" - $\F$-operators.
\begin{prop}
There exists operator $\RR_{1}$ which is the solution of
defining equations
\begin{equation} \RR_{1}
\mathrm{L}_{1}(u_1,u_2)\mathrm{L}_{2}(v_1,v_2)=
\mathrm{L}_{1}(v_1,u_2)\mathrm{L}_{2}(u_1,v_2)\RR_{1}
\label{F1}
\end{equation}
$$
\RR_1 = \RR_1(u_1|v_1,v_2)\ ;\ \RR_1(u_1|v_1,v_2) =
\RR_1(u_1+\lambda|v_1+\lambda,v_2+\lambda)
$$
and these requirements fix the operator $\RR_1$ up to
overall normalization constant
$$ \RR_1(u_1|v_1,v_2) \sim
\frac{\Gamma(z_{21}\dd_2+u_1-v_2)}{\Gamma(z_{21}\dd_2+v_1-v_2)}\
;\ z_{21} = z_2-z_1
$$
\end{prop}
\begin{prop}
There exists operator $\RR_{2}$ which is the solution of
defining equations
\begin{equation}
\RR_{2} \mathrm{L}_{1}(u_1,u_2)\mathrm{L}_{2}(v_1,v_2)=
\mathrm{L}_{1}(u_1,v_2)\mathrm{L}_{2}(v_1,u_2) \RR_{2}
\label{F2}
\end{equation}
$$
\RR_2 = \RR_2(u_1,u_2|v_2)\ ;\ \RR_2(u_1,u_2|v_2) =
\RR_2(u_1+\lambda,u_2+\lambda|v_2+\lambda)
$$
and these requirements fix the operator $\RR_2$ up to
overall normalization constant
$$ \RR_2(u_1,u_2|v_2) \sim
\frac{\Gamma(z_{12}\dd_1+u_1-v_2)}{\Gamma(z_{12}\dd_1+u_1-u_2)}\
;\ z_{12} = z_1-z_2
$$
\end{prop}
\begin{prop}
The operator $\check{\R}$ can be factorized in a following
way
\begin{equation}
\check{\R}(u_{1},u_2|v_{1},v_2)=
\RR_{1}(u_1|v_1,u_2)\RR_{2}(u_1,u_2|v_{2}) =
\RR_{2}(v_1,u_2|v_2)\RR_{1}(u_1|v_1,v_{2}) \label{Rfact}
\end{equation}
\end{prop}
Note that the $\F$-operators change the spins of
$s\ell(2)$-representations \be \F_1 : V_{\ell_1}\otimes
V_{\ell_2} \to V_{\ell_1-\xi_1}\otimes V_{\ell_2+\xi_1} \
;\ \xi_1 = \frac{u_1-v_1}{2} \label{slF1} \ee \be \F_2 :
V_{\ell_1}\otimes V_{\ell_2} \to V_{\ell_1+\xi_2}\otimes
V_{\ell_2-\xi_2} \ ;\ \xi_2 = \frac{u_2-v_2}{2}\label{slF2}
\ee but the general R-matrix $\R_{\ell_1 \ell_2}(u) =
\P_{12}\F_{1}(u_1|v_1,u_2)\F_{2}(u_1,u_2|v_{2})$ appears
automatically $s\ell(2)$-invariant.

The factorization of the $\check{\R}$-operator can be
proven using the simple pictures. The operator $\check{\R}$
interchanges all parameters in the product of two
$\mathrm{L}$-operators. The operator $\F_{2}$ interchanges
the parameters $u_2$ and $v_2$ only and the operator
$\F_{1}$ interchanges the parameters $u_1$ and $v_1$. Using
the operator $\F_{1}\F_{2}$ it is possible to interchange
parameters $u_1,v_1$ and $u_2,v_2$ in two steps so that we
obtain the first equality in~(\ref{Rfact}) as the condition
of commutativity for the diagram

\vspace{5mm} \unitlength 0.8mm \linethickness{0.4pt}
\begin{picture}(130.00,35.00)
\put(50.00,10.00){\vector(1,1){15.00}}
\put(33.00,20.00){\makebox(0,0)[cc]{$\F_{2}(u_{1},u_2|v_2)$}}
\put(85.00,25.00){\vector(1,-1){15.00}}
\put(117.00,20.00){\makebox(0,0)[cc]{$\F_{1}(u_{1}|v_{1},u_2)$}}
\put(25.00,10.00){\makebox(0,0)[cc]{$\mathrm{L}_1(u_1,u_2)
\mathrm{L}_2(v_1,v_2)$}}
\put(75.00,30.00){\makebox(0,0)[cc]{$\mathrm{L}_1(u_1,v_2)
\mathrm{L}_2(v_1,u_2)$}}
\put(125.00,10.00){\makebox(0,0)[cc]{$\mathrm{L}_1(v_1,v_2)
\mathrm{L}_2(u_1,u_2)$}}
\put(75.00,5.00){\makebox(0,0)[cc]{$\check{\R}(u_{1},u_2|v_1,v_{2})$}}
\put(55.00,10.00){\vector(1,0){40.00}}
\end{picture}
\vspace{5mm}

It is possible to exchange the parameters in the opposite
order so that the second equality in~(\ref{Rfact}) is the
condition of commutativity for the diagram

\vspace{5mm} \unitlength 0.8mm \linethickness{0.4pt}
\begin{picture}(130.00,35.00)
\put(50.00,10.00){\vector(1,1){15.00}}
\put(33.00,20.00){\makebox(0,0)[cc]{$\F_{1}(u_{1}|v_1,v_2)$}}
\put(85.00,25.00){\vector(1,-1){15.00}}
\put(117.00,20.00){\makebox(0,0)[cc]{$\F_{2}(v_{1},u_{2}|v_2)$}}
\put(25.00,10.00){\makebox(0,0)[cc]{$\mathrm{L}_1(u_1,u_2)
\mathrm{L}_2(v_1,v_2)$}}
\put(75.00,30.00){\makebox(0,0)[cc]{$\mathrm{L}_1(v_1,u_2)
\mathrm{L}_2(u_1,v_2)$}}
\put(125.00,10.00){\makebox(0,0)[cc]{$\mathrm{L}_1(v_1,v_2)
\mathrm{L}_2(u_1,u_2)$}}
\put(75.00,5.00){\makebox(0,0)[cc]{$\check{\R}(u_{1},u_2;|v_1,v_{2})$}}
\put(55.00,10.00){\vector(1,0){40.00}}
\end{picture}
\vspace{5mm}

The defining system of equations for the $\F$-operator is
the system of differential equations of the second order.
It can be reduced to the simpler system of equations of the
first order which clearly shows the $s\ell(2)$-covariance
of the $\F$-operator.

\begin{lem}
The system of equations~(\ref{F1}) for the operator $\F_1$
is equivalent to the system \be \F_{1}\cdot
\left[\mathrm{L}_1(u_1,u_2)+ \mathrm{L}_2(v_1,v_2)\right] =
\left[ \mathrm{L}_1(v_1,u_2)+
\mathrm{L}_2(u_1,v_2)\right]\cdot \F_1 \ ;\ \F_{1}\cdot z_1
= z_1 \cdot \F_1\label{F1def} \ee
\end{lem}

\begin{lem}
The system of equations~(\ref{F2}) for the operator $\F_2$
is equivalent to the system \be \F_{2}\cdot \left[
\mathrm{L}_1(u_1,u_2)+ \mathrm{L}_2(v_1,v_2)\right] =
\left[ \mathrm{L}_1(u_1,v_2)+
\mathrm{L}_2(v_1,u_2)\right]\cdot \F_2\ ;\ \F_{2}\cdot z_2
= z_2 \cdot \F_2 \label{F2def} \ee
\end{lem}
Note that the relations in the first place are simply the
rules of commutation of $\F$-operators with
$s\ell(2)$-generators written in a compact form. In
explicit notations it is exactly the
relations~(\ref{slF1}),~(\ref{slF2}). The
$s\ell(2)$-invariance of $\R$-matrix follows directly from
the properties of $\F$-operators.\\
{\bf Proof of the Lemma 2 and the Proposition 2.} As an
example we shall consider the operator $\F_2$ and all
calculations for the operator $\F_1$ are very similar. We
are going to prove that the defining equation~(\ref{F2}) is
equivalent to the system~(\ref{F2def}) and derive the
explicit formula for the operator $\F_2$. First we show
that the system~(\ref{F2def}) is the direct consequence of
the eq.~(\ref{F2}). Let us make the shift $u_1\to
u_1+\lambda\ ,\ u_2\to u_2+\lambda\ ,\ v_1\to v_1+\mu\ ,\
v_2\to v_2+\lambda$ in the defining equation ~(\ref{F2}).
The $\F$-operator is invariant under this shift,
$\mathrm{L}$-operators transform simply
$$
\mathrm{L}_1\to \mathrm{L}_1 + \lambda\cdot\II\ ,\
\mathrm{L}_2\to \mathrm{L}_2 + \lambda\cdot\II
+(\mu-\lambda)\left(\begin{array}{cc}
1 & 0 \\
z_2 & 0
   \end{array}\right)
$$
so that we derive as consequence of eq.~(\ref{F2})
$$
\F_2\cdot \left[\lambda\cdot
[\mathrm{L}_1(u_1,u_2)+\mathrm{L}_2(v_1,v_2)]+
(\mu-\lambda)\mathrm{L}_1(u_1,u_2)\left(\begin{array}{cc}
1 & 0 \\
z_2 & 0
   \end{array}\right)+
   \lambda(\mu-\lambda)\left(\begin{array}{cc}
1 & 0 \\
z_2 & 0
   \end{array}\right)\right] =
$$
$$
= \left[\lambda\cdot
[\mathrm{L}_1(u_1,v_2)+\mathrm{L}_2(v_1,u_2)]+
(\mu-\lambda)\mathrm{L}_1(u_1,v_2)\left(\begin{array}{cc}
1 & 0 \\
z_2 & 0
   \end{array}\right)+
   \lambda(\mu-\lambda)\left(\begin{array}{cc}
1 & 0 \\
z_2 & 0
   \end{array}\right)\right]\cdot\F_2
$$
The parameters $\lambda$ and $\mu$ are arbitrary so that
we obtain the equations
$$
\F_{2}(u_1,u_2|v_2)\cdot \left[ \mathrm{L}_1(u_1,u_2)+
\mathrm{L}_2(v_1,v_2)\right] = \left[
\mathrm{L}_1(u_1,v_2)+ \mathrm{L}_2(v_1,u_2)\right]\cdot
\F_2 (u_1,u_2|v_2)\ ;\  \F_{2}\cdot z_2 = z_2 \cdot \F_2\
$$
$$
\F_{2}\cdot\mathrm{L}_1(u_1,u_2) \left(\begin{array}{cc}
1\\
z_{2}
   \end{array}\right) =
   \mathrm{L}_1(u_1,v_2) \left(\begin{array}{cc}
1\\
z_{2}
   \end{array}\right)\cdot\F_2
$$
It is easy to prove that the last equation is the direct
consequence of the first and the second ones. Next we show
that from the systems of equations~(\ref{F2def}) follows
eq.~(\ref{F2}). This will be almost trivial if we rewrite
the equations~(\ref{F2def}) and~(\ref{F2}) in equivalent
form using the $s\ell(2)$-invariance of the
$\mathrm{L}$-operator and the commutativity of $\F_2$ and
$z_2$. We start from the equation~(\ref{F2}) for the
operator $\F_2$ and use the factorized
expression~(\ref{factor}) for the $\mathrm{L}_2$-operator
$$
\F_{2}\cdot \mathrm{L}_1(u_1,u_2) \mathbf{M}
\left(\begin{array}{cc}
v_1-1& -\dd_{2} \\
0& v_2
   \end{array}\right)
\mathbf{M}^{-1}=
   \mathrm{L}_1(u_1,v_2)\mathbf{M} \left(\begin{array}{cc}
v_1-1& -\dd_{2} \\
0& u_2
   \end{array}\right)
\mathbf{M}^{-1}\cdot\F_{2}\ ;\ \mathbf{M}
=\left(\begin{array}{cc}
1 & 0 \\
z_{2}& 1
   \end{array}\right)
$$
Next we perform the similarity transformation
$\mathbf{M}^{-1}\cdots\mathbf{M}$ of this matrix equation
$$
\F_{2}\cdot\mathbf{M}^{-1}\mathrm{L}_1(u_1,u_2)\mathbf{M}
\cdot\left(\begin{array}{cc}
v_1-1& -\dd_{2} \\
0& v_2
   \end{array}\right)
= \mathbf{M}^{-1}\mathrm{L}_1(u_1,v_2)\mathbf{M}
\cdot\left(\begin{array}{cc}
v_1-1& -\dd_{2} \\
0& u_2
   \end{array}\right)\cdot\F_{2}
$$
The $s\ell(2)$-invariance of $\mathrm{L}$-operator allows
to rewrite the matrix $\mathbf{M}^{-1}\mathrm{L}\mathbf{M}$
in the form~(\ref{sl2})
$$
\mathbf{M}^{-1}\mathrm{L}_1(u_1,u_2)\mathbf{M} =
\mathrm{e}^{-z_2\dd_{1}}
   \mathrm{L}_1(u_1,u_2)\mathrm{e}^{+z_2\dd_{1}}
$$
so that we have the simple equation for the operator
$\mathbf{r}\equiv\mathrm{e}^{z_2\dd_{1}}\cdot
\F_{2}\cdot\mathrm{e}^{-z_2\dd_{1}}$~(note that
$\mathbf{r}z_2 = z_2\mathbf{r}$) \be
\mathbf{r}\cdot\mathrm{L}_1(u_1,u_2)
 \left(\begin{array}{cc}
v_1-1& -\dd_{2}+\dd_1 \\
0& v_2
   \end{array}\right)
= \mathrm{L}_1(u_1,v_2)\left(\begin{array}{cc}
v_1-1& -\dd_{2}+\dd_1 \\
0& u_2
   \end{array}\right)\cdot\mathbf{r}
\label{def} \ee This system of equations is equivalent to
the system~(\ref{F2}) written in terms of $\mathbf{r}$. To
derive the system of equations which is equivalent to the
system~(\ref{F2def}) written in terms of $\mathbf{r}$ we
repeat the same trick with the shift of parameters and
obtain the equations \be \mathbf{r} \left[
\mathrm{L}_1(u_1,u_2)+\left(\begin{array}{cc}
v_1-1& -\dd_{2}+\dd_1 \\
0& v_2
   \end{array}\right)\right] =
\left[\mathrm{L}_1(u_1,v_2)+ \left(\begin{array}{cc}
v_1-1& -\dd_{2}+\dd_1 \\
0& u_2
   \end{array}\right)\right]\mathbf{r}
\label{Def}\ee
$$
\mathbf{r}\cdot\mathrm{L}_1(u_1,u_2)
\left(\begin{array}{cc} 1 \\ 0
   \end{array}\right) =
\mathrm{L}_1(u_1,v_2) \left(\begin{array}{cc} 1  \\0
   \end{array}\right)\cdot\mathbf{r}
$$
The second equation is the evident consequence of the first
one. We must prove that the system~(\ref{def}) is
equivalent to the system~(\ref{Def}). First step we
factorize in~(\ref{def}) the matrix $diag(v_1-1\ ;\ 1)$
from the right so that the parameter $v_1$ disappears from
the equation
$$
\mathbf{r}\cdot\mathrm{L}_1(u_1,u_2)
 \left(\begin{array}{cc}
1& -\dd_{2}+\dd_1 \\
0& v_2
   \end{array}\right)=
   \mathrm{L}_1(u_1,v_2)
\left(\begin{array}{cc}
1& -\dd_{2}+\dd_1 \\
0& u_2
   \end{array}\right)\cdot\mathbf{r}
$$
It is easy to see that there are two new equations in
comparison with~(\ref{Def})
\be
\mathbf{r}\cdot\mathrm{L}_1(u_1,u_2)
 \left(\begin{array}{cc}
\dd_1 \\
v_2
   \end{array}\right)=
   \mathrm{L}_1(u_1,v_2)
\left(\begin{array}{cc}
\dd_1 \\
u_2
   \end{array}\right)\cdot\mathbf{r}
\label{new}\ee In explicit form the defining
system~(\ref{Def}) contains three equations \be \mathbf{r}
z_1\dd_1 = z_1\dd_1 \mathbf{r} \ ;\ \mathbf{r} \dd_2 =
\dd_2 \mathbf{r} \ ;\ \mathbf{r} z_1(z_1\dd_1+u_1-u_2) =
z_1(z_1\dd_1+u_1-v_2)\mathbf{r} \label{Deff}\ee The
"down"-equation in~(\ref{new}) is the simple consequence of
equation $\mathbf{r}\cdot z_1\dd_1 =
z_1\dd_1\cdot\mathbf{r}$ from the system~(\ref{Deff}) and
the "up"-equation
$$
\mathbf{r}(z_1\dd_1+u_1-v_2)\dd_1 =
(z_1\dd_1+u_1-u_2)\dd_1\mathbf{r}
$$
is equivalent to the last equation from the
system~(\ref{Deff}). Indeed we have
$$
\mathbf{r} z_1(z_1\dd_1+u_1-u_2) =
z_1(z_1\dd_1+u_1-v_2)\mathbf{r} \rightarrow \mathbf{r}
z_1(z_1\dd_1+u_1-u_2)\dd_1 =
z_1(z_1\dd_1+u_1-v_2)\mathbf{r}\dd_1 \rightarrow
$$
$$
\rightarrow z_1(z_1\dd_1+u_1-u_2)\dd_1\mathbf{r} =
z_1\mathbf{}(z_1\dd_1+u_1-v_2)\dd_1 \rightarrow
(z_1\dd_1+u_1-u_2)\dd_1\mathbf{r} =
\mathbf{r}(z_1\dd_1+u_1-v_2)\dd_1
$$
We have proved the equivalence of the systems~(\ref{def})
and~(\ref{Def}) and therefore the equivalence of the
systems~(\ref{F2def}) and~(\ref{F2}). It remains to find
the solution. First we solve the equations~(\ref{Deff}) for
operator $\mathbf{r}$. The solution of equation
$\mathbf{r}z_2 = z_2\mathbf{r}$ and the first two equations
from the system~(\ref{Deff}) is $\mathbf{r} =
\mathbf{r}[z_1\dd_1]$. Then the last equation leads to the
recurrence relation which has the simple solution
$$
\mathbf{r}[z_1\dd_1+1](z_1\dd_1+u_1-u_2)=
\mathbf{r}[z_1\dd_1](z_1\dd_1+u_1-v_2)\Longrightarrow
\mathbf{r}[z_1\dd_1] \sim
\frac{\Gamma(z_1\dd_1+u_1-v_2)}{\Gamma(z_1\dd_1+u_1-u_2)}
$$
Finally we derive the expression for the operator $\F_{2}$
from the Proposition
$$
\F_{2}=
\mathrm{e}^{-z_2\dd_{1}}\mathbf{r}\textrm{e}^{z_2\dd_{1}}
\sim
\frac{\Gamma(z_{12}\dd_1+u_1-v_2)}{\Gamma(z_{12}\dd_1+u_1-u_2)}.
$$

\section{The general $sl(3)$-invariant R-matrix}
\setcounter{equation}{0}

In the previous section we have constructed the general
solution of the Yang-Baxter equation with symmetry algebra
$s\ell(2)$ and proved that it has simple factorized
structure. Next we shall consider the first nontrivial
example of the symmetry algebra of the rank two. In fact,
we repeat step by step all calculations from the previous
section and show that the general solution of the
Yang-Baxter equation with the symmetry algebra $s\ell(3)$
has the very similar structure.

\subsection{$sl(3)$ lowest weight modules}
\setcounter{equation}{0}

The algebra $s\ell(3)$ has eight generators
$\mathbf{T}_{ab}\ ,\ a,b = 1 , 2 , 3$ with condition
$\mathbf{T}_{11}+\mathbf{T}_{22}+\mathbf{T}_{33} = 0$. The
commutation relations have the standard form~\cite{Z,dS}
$$
[\mathbf{T}_{ab},\mathbf{T}_{cd}] = \delta_{cb}
\mathbf{T}_{ad}-\delta_{ad}\mathbf{T}_{cb}
$$
We shall use the following generators of the Cartan
subalgebra
$$ \mathbf{H}_1 = \mathbf{T}_{11}-\mathbf{T}_{22}\
;\ \mathbf{H}_2 = \mathbf{T}_{22}-\mathbf{T}_{33}\ ;\
\mathbf{T}_{11} = \frac{2}{3}\mathbf{H}_1 +
\frac{1}{3}\mathbf{H}_2\ ;\ \mathbf{T}_{22} =
\frac{1}{3}\mathbf{H}_2 - \frac{1}{3}\mathbf{H}_1 \ ;\
\mathbf{T}_{33} = -\frac{1}{3}\mathbf{H}_1 -
\frac{2}{3}\mathbf{H}_2
$$
There are two central elements $\mathbf{C}_2 = \sum_{ab}
\mathbf{T}_{ab}\mathbf{T}_{ba}$ and $\mathbf{C}_3 =
\sum_{abc} \mathbf{T}_{ab}\mathbf{T}_{bc}\mathbf{T}_{ca}$.
The Verma module is the generic lowest weight
$s\ell(3)$-module $\mathbf{V}_{\Lambda}\ ;\ \Lambda = (m ,
n) $. As a linear space $\mathbf{V}_{\Lambda}$ is obtained
by application of operators $\mathbf{T}_{12} ,
\mathbf{T}_{13} , \mathbf{T}_{23}$ to the lowest weight
vector $\mathbf{a}_0$
$$
\mathbf{T}_{21}\mathbf{a}_0
=\mathbf{T}_{31}\mathbf{a}_0 = \mathbf{T}_{32}\mathbf{a}_0
 = 0 \ ;\ \mathbf{H}_1\mathbf{a}_0 = -m\cdot\mathbf{a}_0
 \ ;\  \mathbf{H}_2\mathbf{a}_0 = -n\cdot\mathbf{a}_0
$$
We shall use the representation $\mathrm{V}_{\Lambda}$ of
$s\ell(3)$ in the infinite-dimensional space $\C[x,y,z]$ of
polynomials in variables $x,y,z$ and lowest weight vector
$a_0 = 1$~\cite{Z,Shif}. The action of $s\ell(3)$ in
$\mathrm{V}_{\Lambda}$ is given by the first-order
differential operators. Lowering (decreasing the polynomial
degree) operators have the form \be \mathrm{T}_{21} =
\partial_{x} \ ;\ \mathrm{T}_{31} =
\partial_{y} \ ;\ \mathrm{T}_{32}
= \partial_{z}-x\partial_{y} \label{gen3-} \ee and generate
the following global transformations
$$
\mathrm{e}^{\lambda \mathrm{T}_{21}}\Phi(x,y,z)=
\Phi(x+\lambda,y,z) \ ;\
\mathrm{e}^{\lambda\mathrm{T}_{31}}\Phi(x,y,z)=
\Phi\left(x,y+\lambda,z\right)
$$
$$
\mathrm{e}^{\lambda\mathrm{T}_{32}}\Phi(x,y,z)= \Phi\left(x
, y-\lambda x , z+\lambda\right).
$$
Rising(increasing the polynomial degree) operators
$$
\mathrm{T}_{12}= -x^2\partial_{x} - xy\partial_{y} +x
z\partial_{z} + y\partial z + n x \ ,\ \mathrm{T}_{23}=
-z^2\partial_z -y\partial_{x} + m z
$$
\be\mathrm{T}_{13}= -y^2\partial_{y} -x y\partial_{x} -
z(y+xz)\partial_{z}+(m+n)y + m x z \label{gen3+} \ee
generate the global transformations
$$
\mathrm{e}^{\lambda\mathrm{T}_{12}}\Phi(x,y,z)=
\left[1+\lambda x\right]^{n}\Phi\left(\frac{x}{1+\lambda
x},\frac{y}{1+\lambda x}, z+\lambda(y+xz)\right)
$$
$$
\mathrm{e}^{\lambda\mathrm{T}_{23}}\Phi(x,y,z)=
\left[1+\lambda z\right]^{m} \Phi\left(x-\lambda y , y
,\frac{z}{1+\lambda z}\right) ,
$$
$$
\mathrm{e}^{\lambda\mathrm{T}_{13}}\Phi(x,y,z) =
\left[1+\lambda
y\right]^{n}\left[1+\lambda(y+xz)\right]^{m}
\Phi\left(\frac{x}{1+\lambda y} , \frac{y}{1+\lambda y},
\frac{z}{1+\lambda(y+xz)}\right) .
$$
Two remaining elements of the Cartan subalgebra:
$$
\mathrm{H}_1 = 2x\partial_x +y\partial_y-z\partial_z -n \
;\ \mathrm{H}_2 = 2z\partial_z +y\partial_y-x\partial_x - m
$$
generate the transformations:
$$
\lambda^{\mathrm{H}_1}\Phi(x,y,z) = \lambda^{-n}
\Phi\left(\lambda^2 x , \lambda y ,
\frac{z}{\lambda}\right)\ ;\
\lambda^{\mathrm{H}_2}\Phi(x,y,z) = \lambda^{-m}
\Phi\left(\frac{z}{\lambda} , \lambda y , \lambda^2
z\right)
$$
Using these formulae it is possible to derive the closed
expression for the generating functions of the basis
vectors
$$
\mathrm{e}^{\mu\mathrm{T}_{12}}
\mathrm{e}^{\nu\mathrm{T}_{23}}
\mathrm{e}^{\lambda\mathrm{T}_{13}}\cdot 1 = \left[1+\mu
x+\lambda y\right]^{n}\cdot\left[1+\nu
z+(\lambda+\mu\nu)(y+xz)\right]^{m}
$$
The power expansion of these generating functions in
$\mu,\nu ,\lambda$ gives the elements of the basis. It is
evident that for the generic $m\neq M\ ,\ m\neq N \ ;\ M,N
\in\{0,1,2,3\cdots\}$ the module $\mathrm{V}_{\Lambda}$ is
an irreducible lowest weight $s\ell(3)$-module isomorphic
to $\mathbf{V}_{\Lambda}$ but for the special values of the
spin $m = M\ ,\ n = N$ there exists the finite dimensional
invariant subspace $\mathrm{V}_{M,N}$ with dimension $\dim
\mathrm{V}_{M,N} =
\frac{1}{2}(M+1)(N+1)(M+N+2)$~\cite{Z,dS}. For generic $M,N
\ne 0$ the space $\mathrm{V}_{M,N}$ is invariant subspace
in the space of polynomials in three variables $\C[x,y,z]$.
In the case $N = 0$ we have
$$
\mathrm{e}^{\mu\mathrm{T}_{12}}
\mathrm{e}^{\nu\mathrm{T}_{23}}
\mathrm{e}^{\lambda\mathrm{T}_{13}}\cdot 1 = \left[1+\nu
z+(\lambda+\mu\nu)(y+xz)\right]^{M}
$$
so that one obtains the invariant subspace
$\mathrm{V}_{M,0}$ with dimension $\dim \mathrm{V}_{M,0} =
\frac{1}{2}(M+1)(M+2)$ in the space of polynomials in two
variables $\C[y+xz,z]$. In the case $M = 0$ we have
$$
\mathrm{e}^{\mu\mathrm{T}_{12}}
\mathrm{e}^{\nu\mathrm{T}_{23}}
\mathrm{e}^{\lambda\mathrm{T}_{13}}\cdot 1 = \left[1+\mu
x+\lambda y\right]^{N}
$$
so that one obtains the invariant subspace
$\mathrm{V}_{0,N}$ with dimension $\dim \mathrm{V}_{0,N} =
\frac{1}{2}(N+1)(N+2)$ in the space of polynomials in two
variables $\C[x , y]$. We shall use the three-dimensional
representation $\mathrm{V}_{(1,0)}\sim\C^3$.In the basis
$$
\mathbf{e}_1 = \mathrm{T}_{13}\cdot 1 = y+x z\ ,\
\mathbf{e}_2 = \mathrm{T}_{23}\cdot 1 = z\ ,\ \mathbf{e}_3
= 1
$$
the $s\ell(3)$-generators take the form
\be
\mathbf{t}_{31}
= \left (\begin{array}{ccc}
0 & 0 & 0  \\
0 & 0 & 0  \\
1 & 0 & 0
\end{array} \right )
\ ;\  \mathbf{t}_{21} = \left (\begin{array}{ccc}
0 & 0 & 0  \\
1 & 0 & 0  \\
0 & 0 & 0
\end{array} \right )
\ ;\  \mathbf{t}_{32} = \left (\begin{array}{ccc}
0 & 0 & 0  \\
0 & 0 & 0  \\
0 & 1 & 0
\end{array} \right )\ ;\ \mathbf{h}_1 = \left (\begin{array}{ccc}
1 & 0 & 0  \\
0 & -1 & 0  \\
0 & 0 & 0
\end{array} \right )
\label{fund1} \ee
$$
\mathbf{t}_{13} = \left (\begin{array}{ccc}
0 & 0 & 1  \\
0 & 0 & 0  \\
0 & 0 & 0
\end{array} \right )
\ ;\  \mathbf{t}_{23} = \left (\begin{array}{ccc}
0 & 0 & 0  \\
0 & 0 & 1  \\
0 & 0 & 0
\end{array} \right )
\ ;\  \mathbf{t}_{12} = \left (\begin{array}{ccc}
0 & 1 & 0  \\
0 & 0 & 0  \\
0 & 0 & 0
\end{array} \right )\ ;\  \mathbf{h}_2 = \left (\begin{array}{ccc}
0 & 0 & 0  \\
0 & 1 & 0  \\
0 & 0 & -1
\end{array} \right )
$$
There exists the second three-dimensional representation
$\mathrm{V}_{(0,1)}\sim\C^3$. In the basis
$$
\mathbf{e}_1 = -1\ ;\ \mathbf{e}_2 = \mathrm{T}_{12}\cdot 1
= x\ ,\ \mathbf{e}_3 = \mathrm{T}_{13}\cdot 1 = y $$ the
$s\ell(3)$-generators take the similar form but
$\mathbf{t}_{ik}\to -\mathbf{t}_{ki}\ ;\ \mathbf{h}_1\to
-\mathbf{h}_1\ ;\ \mathbf{h}_2\to -\mathbf{h}_2$.

\subsection{Yang-Baxter equation and Lax operator}
\setcounter{equation}{0}

The Yang-Baxter equation is the following three term
relation
$$
\R_{\Lambda_1\Lambda_2}(u-v)\R_{\Lambda_1\Lambda_3}(u)
\R_{\Lambda_2\Lambda_3}(v)=
\R_{\Lambda_2\Lambda_3}(v)\R_{\Lambda_1\Lambda_3}(u)
\R_{\Lambda_1\Lambda_2}(u-v)
$$
for the operators $\R_{\Lambda_i\Lambda_j}(u):
V_{\Lambda_i}\otimes V_{\Lambda_j}\to V_{\Lambda_i}\otimes
V_{\Lambda_j}$. As in the $s\ell(2)$-case we start from the
simplest solutions of Yang-Baxter equation and derive the
defining equation for the general $\R$-operator. First we
put $\Lambda_1 = \Lambda_2 = \Lambda_3 = (1,0)$ in
Yang-Baxter equation and consider the restriction to the
invariant subspace $\C^3\otimes\C^3\otimes\C^3$. We obtain
the equation \be \R_{1 2}(u-v)\R_{1 3}(u) \R_{2 3}(v)=
\R_{2 3}(v)\R_{1 3}(u) \R_{1 2}(u-v) \ee where the operator
$\R_{12}(u)$ acts on the first and second copy of $\C^3$ in
the tensor product $\C^3\otimes\C^3\otimes\C^3$ and
similarly for the other $\R$-operators. The solution is
well known~\cite{KS,KR}
$$
\R_{12}(u) = u + \mathrm{P}_{12}
$$
where $\mathrm{P}_{12}$ is the permutation operator in
$\C^3\otimes\C^3$. Secondly we choose $\Lambda_1 =
\Lambda_2 = (1,0)\ ;\ \Lambda_3 = \Lambda = (m,n)$ and
consider the restriction to the invariant subspace
$\C^3\otimes\C^3\otimes\mathrm{V}_{\Lambda}$. The
restriction of the operator $\R_{\Lambda_1 \Lambda}(u)$ to
the space $\C^3\otimes\mathrm{V}_{\Lambda}$ coincides up to
normalization and shift of spectral parameter with the
Lax-operator
$$
\mathrm{L}(u): \C^3\otimes\mathrm{V}_{\Lambda}\to
\C^3\otimes\mathrm{V}_{\Lambda}
$$
and the Yang-Baxter equation coincides with the fundamental
commutation relations for the
Lax-operator~\cite{KS,KR,M,MNO}
$$
\R_{12}(u-v)\mathrm{L}^{(1)}(u)\mathrm{L}^{(2)}(v) =
\mathrm{L}^{(2)}(v)\mathrm{L}^{(1)}(u)\R_{12}(u-v)
$$
where $\mathrm{L}^{(1)}(u)$ is the operator which acts
nontrivially on the first copy of $\C^3$ and
$\mathrm{V}_{\Lambda}$ in the tensor product
$\C^3\otimes\C^\otimes\mathrm{V}_{\Lambda}$ and
$\mathrm{L}^{(2)}(u)$ is the operator which acts
nontrivially on the second copy of $\C^3$ and
$\mathrm{V}_{\Lambda}$. The solution coincides up to
additive constant with the Casimir operator $\mathbf{C}_2$
for the representation $\C^3\otimes\mathrm{V}_{\Lambda}$
$$
\mathrm{L}(u) \equiv  u  + \frac{1}{3}\cdot
\left(2\mathbf{h}_1 +
\mathbf{h}_2\right)\otimes\mathrm{H}_1 + \frac{1}{3}\cdot
\left(\mathbf{h}_1 +
2\mathbf{h}_2\right)\otimes\mathrm{H}_2 + \sum_{i\ne
k}\mathbf{t}_{ik}\otimes\mathrm{T}_{ki}
$$
where $\mathbf{h}_1 , \mathbf{h}_2 , \mathbf{t}_{ik}$ are
$s\ell(3)$-generators in the fundamental representation and
$\mathrm{H}_1 , \mathrm{H}_2 , \mathrm{T}_{ik}$ are
generators in the generic representation. The algebra
$s\ell(3)$ has two three-dimensional representations --
$\mathrm{V}_{(1,0)}$ and $\mathrm{V}_{(0,1)}$ so that there
exists the second Lax-operator
$$
\bar{\mathrm{L}}(u): \C^3\otimes\mathrm{V}_{\Lambda}\to
\C^3\otimes\mathrm{V}_{\Lambda}\ ;\
\R_{12}(u-v)\bar{\mathrm{L}}^{(1)}(u)\bar{\mathrm{L}}^{(2)}(v)
=
\bar{\mathrm{L}}^{(2)}(v)\bar{\mathrm{L}}^{(1)}(u)\R_{12}(u-v)
$$
The explicit expression for the second Lax-operator is the
same but now $\mathbf{t}_{ik}\to -\mathbf{t}_{ki}\ ;\
\mathbf{h}_1\to -\mathbf{h}_1\ ;\ \mathbf{h}_2\to
-\mathbf{h}_2$. We shall use the Lax-operator
$\mathrm{L}(u)$
$$
\mathrm{L}(u) = \left(\begin{array}{ccc}
\frac{2}{3}\mathrm{H}_1+\frac{1}{3}\mathrm{H}_2+u & \mathrm{T}_{21} & \mathrm{T}_{31}\\
\mathrm{T}_{12} & \frac{1}{3}\mathrm{H}_2-\frac{1}{3}\mathrm{H}_1+u & \mathrm{T}_{32}\\
\mathrm{T}_{13} & \mathrm{T}_{23} &
-\frac{1}{3}\mathrm{H}_1-\frac{2}{3}\mathrm{H}_2+u
   \end{array}\right)
$$
in defining equation for the general $\R$-operator. The
Lax-operator $\mathrm{L}(u)$ depends on three parameters $u
, m , n$. We shall use the parametrization
$$u_1 = u-2-\frac{m+2n}{3}\ ,\ u_2 = u-1+\frac{n-m}{3}
\ ,\ u_3 = u+\frac{n+2m}{3}\ ;\ m = u_3-u_2-1\ ,\ n =
u_2-u_1-1
$$
and show this parameters explicitly. The Lax operator
$\mathrm{L}(u_1,u_2,u_3)$ in the functional representation
$\mathrm{V}_{\Lambda}$ has the form \be
\mathrm{L}(u_1,u_2,u_3) =\left(\begin{array}{ccc}
x\dd_x+y\dd_y+u_1+2 & \dd_x & \dd_y\\
\mathrm{L}_{21}&
-x\dd_x+z\dd_z+u_2+1 & \dd_z-x\dd_y\\
\mathrm{L}_{31} & \mathrm{L}_{32} & -y\dd_y-z\dd_z+u_3
   \end{array}\right)
\label{Lax3} \ee
$$
\mathrm{L}_{21} =-x^2\dd_x
-xy\dd_y+(xz+y)\dd_z+(u_2-u_1-1)x \ ;\ \mathrm{L}_{32} =
-y\dd_x-z^2\dd_z+(u_3-u_2-1)z
$$
$$
\mathrm{L}_{31} =
-xy\dd_x-y^2\dd_y-z(xz+y)\dd_z+(u_3-u_2-1)xz+(u_3-u_1-2)y
$$
and similar to $s\ell(2)$-case there exists the factorized
representation for the Lax-operator
\begin{equation}
\mathrm{L}(u_1,u_2,u_3) \equiv \left(\begin{array}{ccc}
1 & 0 & 0\\
-x & 1 & 0\\
-y & -z & 1
   \end{array}\right)
   \left(\begin{array}{ccc}
u_1 & \dd_x-z\dd_y & \dd_y \\
0 & u_2 & \dd_z\\
0 & 0 & u_3
   \end{array}\right)
   \left(\begin{array}{ccc}
1 & 0 & 0\\
x & 1 & 0\\
y+xz & z & 1
   \end{array}\right)
\label{factor3}
\end{equation}
The $\mathrm{L}$-operator is $s\ell(3)$-invariant by
construction and as consequence one obtains the useful
equality \be
\mathbf{M}^{-1}\cdot\mathrm{L}(u)\cdot\mathbf{M} =
\S^{-1}\cdot\mathrm{L}(u)\cdot\S\ ;\ \S =
\mathrm{e}^{c\left(\dd_{z}-x\dd_{y}\right)}
\cdot\mathrm{e}^{b\dd_{y}} \cdot\mathrm{e}^{a\dd_{x}} \ ;\
\mathbf{M} = \left(\begin{array}{ccc}
1 & 0 & 0\\
-a & 1 & 0\\
-b & -c & 1
   \end{array}\right)
\label{sl3} \ee Finally we put $\Lambda_1 = \left(1,
0\right)$ in Yang-Baxter equation
$$
\R_{\Lambda_1\Lambda_2}(u-v)\R_{\Lambda_1\Lambda_3}(u)
\R_{\Lambda_2\Lambda_3}(v)=
\R_{\Lambda_2\Lambda_3}(v)\R_{\Lambda_1\Lambda_3}(u)
\R_{\Lambda_1\Lambda_2}(u-v)
$$
change the numeration of the representation spaces
$\Lambda_2 \to \Lambda_1 = (m_1, n_1)\ ;\ \Lambda_3 \to
\Lambda_2 = (m_2, n_2)$ and consider the restriction to the
invariant subspace $\C^3\otimes\mathrm{V}_{\Lambda_1}
\otimes\mathrm{V}_{\Lambda_2}$. In this way one obtains the
defining equation for the $\R$-operator
$$
\mathrm{L}_1(u-v)\mathrm{L}_2(u)
\R_{\Lambda_1\Lambda_2}(v)=
\R_{\Lambda_1\Lambda_2}(v)\mathrm{L}_2(u)\mathrm{L}_1(u-v)
$$
The operator $\mathrm{L}_k$ acts nontrivially on the tensor
product $\C^3\otimes\mathrm{V}_{\Lambda_k}$ which is
isomorphic to $\C^3\otimes\C[x_k,y_k,z_k]$ and the operator
$\R_{\Lambda_1\Lambda_2}(u)$ acts nontrivially on the
tensor product
$\mathrm{V}_{\Lambda_1}\otimes\mathrm{V}_{\Lambda_2}$ which
is isomorphic to $\C[x_1,y_1,z_1]\otimes \C[x_2,y_2,z_2]$.
Note that obtained defining equation is slightly different
from the ones we have used in $s\ell(2)$-case. The defining
equation which is similar to $s\ell(2)$ case is \be
\R^{-1}_{\Lambda_1\Lambda_2}(v-u)
\mathrm{L}_1(u)\mathrm{L}_2(v)=
\mathrm{L}_2(v)\mathrm{L}_1(u)
\R^{-1}_{\Lambda_1\Lambda_2}(v-u) \label{def3} \ee There
exists the well known automorphism~\cite{M,MNO} of the
Yang-Baxter equation $\R_{\Lambda_1\Lambda_2}(u) \to
\R^{-1}_{\Lambda_1\Lambda_2}(-u)$. In the simplest
$s\ell(2)$ case we have $\R_{\ell_1\ell_2}(u) \sim
\R^{-1}_{\ell_1\ell_2}(-u)$ but for the more complicated
algebras the action of this automorphism is nontrivial. To
proceed in close analogy with $s\ell(2)$ case we shall use
the defining equation~(\ref{def3}) so that we derive the
expression for the operator
$\R^{-1}_{\Lambda_1\Lambda_2}(v-u)$.

\subsection{The general R-matrix}
\setcounter{equation}{0}

Now we are going to the solution of the defining equation
for the general R-matrix. It is useful to extract the
operator of permutation
$$
\P_{12}: \C[x_1,y_1,z_1]\otimes \C[x_2,y_2,z_2] \to
\C[x_2,y_2,z_2]\otimes \C[x_1,y_1,z_1]
$$
$$
\P_{12}\Psi\left(x_1,y_1,z_1|x_2,y_2,z_2\right) =
\Psi\left(x_2,y_2,z_2|x_1,y_1,z_1\right)
$$
from the $\R$-operator $\R^{-1}_{\Lambda_1\Lambda_2}(v-u) =
\P_{12}\check{\R}(u;v)$ and solve the defining equation for
the $\check{\R}$-operator. The main defining equation for
the $\check{\R}$-operator is
\begin{equation}
\check{\R}(u;v)
\mathrm{L}_{1}(u_1,u_2,u_3)\mathrm{L}_{2}(v_1,v_2,v_3)=
\mathrm{L}_{1}(v_1,v_2,v_3)\mathrm{L}_{2}(u_1,u_2,u_3)
\check{\R}(u;v)
\end{equation}
$$
u_1 = u-2-\frac{m_1+2n_1}{3}\ ,\ u_2 =
u-1+\frac{n_1-m_1}{3} \ ,\ u_3 = u+\frac{n_1+2m_1}{3}
$$
$$
v_1 = v-2-\frac{m_2+2n_2}{3}\ ,\ v_2 =
v-1+\frac{n_2-m_2}{3} \ ,\ v_3 = v+\frac{n_2+2m_2}{3}
$$
The operator $\check{\R}$ can be represented as the product
of the simpler "elementary building blocks" -
$\F$-operators.
\begin{prop}
There exists operator $\F_{1}$ which is the solution of the
defining equations
\begin{equation}
\F_{1}
\mathrm{L}_{1}(u_1,u_2,u_3)\mathrm{L}_{2}(v_1,v_2,v_3)=
\mathrm{L}_{1}(v_1,u_2,u_3)\mathrm{L}_{2}(u_1,v_2,v_3)\F_{1}
\label{3F1}
\end{equation}
$$
\F_1 = \F_1(u_1|v_1,v_2,v_3)\ ;\ \F_1(u_1|v_1,v_2,v_3) =
\F_1(u_1+\lambda|v_1+\lambda,v_2+\lambda,v_3+\lambda)
$$
and these requirements fix the operator $\F_1$ up to
overall normalization constant
$$
\F_1 \sim \S_1^{-1}\cdot \frac{\Gamma(x\dd_{x}+u_1-v_2+1)}
{\Gamma(x\dd_{x}+1)}\cdot \mathrm{e}^{-\frac{y}{x}\dd_{z}}
\cdot \frac{\Gamma(y\partial_{y}+u_1-v_3+1)}
{\Gamma(y\partial_{y}+v_1-v_3+1)}\cdot
\mathrm{e}^{\frac{y}{x}\dd_{z}}\cdot
\frac{\Gamma(x\dd_{x}+1)}
{\Gamma(x\dd_{x}+v_1-v_2+1)}\cdot\S_1
$$
$$
\S_1 = \mathrm{e}^{\left(y_1-z_1x_2\right)\dd_{y_2}}\cdot
\mathrm{e}^{z_1\partial_{z_2}}\cdot
\mathrm{e}^{x_1\dd_{x_2}}
$$
$$
x = x_2\ ,\ y = y_2+x_2 z_2\ ,\ z = z_2\ ;\ \dd_{x} =
\dd_{x_2}-z_2\dd_{y_2}\ ,\ \dd_{y} = \dd_{y_2}\ ,\ \dd_{z}
= \dd_{z_2} -x_2\dd_{y_2}
$$
\end{prop}

\begin{prop}
There exists operator $\F_{2}$ which is the solution of the
defining equations
\begin{equation}
\F_{2}
\mathrm{L}_{1}(u_1,u_2,u_3)\mathrm{L}_{2}(v_1,v_2,v_3)=
\mathrm{L}_{1}(u_1,v_2,u_3)\mathrm{L}_{2}(v_1,u_2,v_3)\F_{2}
\label{3F2}
\end{equation}
$$
\F_2 = \F_2(u_1,u_2|v_2,v_3)\ ;\ \F_2(u_1,u_2|v_2,v_3) =
\F_1(u_1+\lambda,u_2+\lambda|v_2+\lambda,v_3+\lambda)
$$
and these requirements fix the operator $\F_2$ up to
overall normalization constant
$$
\F_2 \sim \S_2^{-1}\cdot
\frac{\Gamma(z_2\partial_{z_2}+u_2-v_3+1)}
{\Gamma(z_2\partial_{z_2}+1)}\cdot
\mathrm{e}^{\frac{y_1}{z_{2}}\dd_{x_1}}\cdot
\frac{\Gamma(x_1\partial_{x_1}+u_1-v_2+1)}
{\Gamma(x_1\partial_{x_1}+u_1-u_2+1)}\cdot
\mathrm{e}^{-\frac{y_1}{z_{2}}\dd_{x_1}}\cdot
\frac{\Gamma(z_2\partial_{z_2}+1)}
{\Gamma(z_2\partial_{z_2}+v_2-v_3+1)}\cdot\S_2
$$
$$
\S_2 = \mathrm{e}^{\left(y_2-x_1 z_1\right)\dd_{y_1}}\cdot
\mathrm{e}^{z_1\partial_{z_2}}\cdot
\mathrm{e}^{x_2\dd_{x_1}}
$$

\end{prop}

\begin{prop}
There exists operator $\F_{3}$ which is the solution of the
defining equations
\begin{equation}
\F_{3}
\mathrm{L}_{1}(u_1,u_2,u_3)\mathrm{L}_{2}(v_1,v_2,v_3)=
\mathrm{L}_{1}(u_1,u_2,v_3)\mathrm{L}_{2}(v_1,v_2,u_3)
\F_{3} \label{3F3}
\end{equation}
$$
\F_3 = \F_3(u_1,u_2,u_3|v_3)\ ;\ \F_3(u_1,u_2,u_3|v_3) =
\F_3(u_1+\lambda,u_2+\lambda,u_3+\lambda|v_3+\lambda)
$$
and these requirements fix the operator $\F_3$ up to
overall normalization constant
$$
\F_3 \sim \S_3^{-1}\cdot
\frac{\Gamma(z_1\partial_{z_1}+u_2-v_3+1)}
{\Gamma(z_1\partial_{z_1}+1)}\cdot
\mathrm{e}^{\frac{y_1}{z_{1}}\dd_{x_1}}\cdot
\frac{\Gamma(y_1\partial_{y_1}+u_1-v_3+1)}
{\Gamma(y_1\partial_{y_1}+u_1-u_3+1)}\cdot
\mathrm{e}^{-\frac{y_1}{z_{1}}\dd_{x_1}}\cdot
\frac{\Gamma(z_1\partial_{z_1}+1)}
{\Gamma(z_1\partial_{z_1}+u_2-u_3+1)} \cdot\S_3
$$
$$
\S_3 = \mathrm{e}^{\left(y_2-x_1 z_2\right)\dd_{y_1}}
\mathrm{e}^{z_2\dd_{z_1}}\cdot \mathrm{e}^{x_2\dd_{x_1}}
$$
\end{prop}

\begin{prop}
The $\check{\R}$-operator can be factorized as follows
\begin{equation}
\check{\R}(u;v)= \F_1(u_1;v_1,u_2,u_3)
\F_2(u_1,u_2;v_2,u_3)\F_3(u_1,u_2,u_3;v_3) \label{Rfact3}
\end{equation}
\end{prop}
There exist six equivalent ways to represent $\check{\R}$
in an factorized form which differ by the order of
$\F$-operators and their parameters. All these expressions
and the proof of the factorization of the
$\check{\R}$-operator can be obtained using the pictures
similar to $s\ell(2)$-case.

The defining systems of equations for the $\F$-operator can
be reduced to the simpler system which clearly shows the
property of $s\ell(3)$-covariance of the $\F$-operator.

\begin{lem}
The defining equation~(\ref{3F1}) for the operator $\F_1$
is equivalent to the system of equations
\begin{equation}
\F_{1} \left[\mathrm{L}_{1}(u_1,u_2,u_3)+
\mathrm{L}_{2}(v_1,v_2,v_3)\right]=
\left[\mathrm{L}_{1}(v_1,u_2,u_3)+
\mathrm{L}_{2}(u_1,v_2,v_3)\right] \F_{1} \label{3F1def}
\end{equation}
$$
\F_{1} x_1 = x_1 \F_{1}\ ,\ \F_{1} y_1 = y_1 \F_{1}\ ,\
\F_{1} z_1 = z_1 \F_{1}
$$
\be\F_1\cdot\left(\dd_{z_2}-x_{21}\dd_{y_2}\right) = \left(
\dd_{z_2}-x_{21}\dd_{y_2}\right)\cdot\F_1 \label{3F1last}
\end{equation}
\end{lem}

\begin{lem}
The defining equation~(\ref{3F2}) for the operator $\F_2$
is equivalent to the system of equations
\begin{equation}
\F_{2} \left[\mathrm{L}_{1}(u_1,u_2,u_3)+
\mathrm{L}_{2}(v_1,v_2,v_3)\right]=
\left[\mathrm{L}_{1}(u_1,v_2,u_3)+
\mathrm{L}_{2}(v_1,u_2,v_3)\right] \F_{2} \label{3F2def}
\end{equation}
$$
\F_{2}\left(y_1+x_1 z_1\right) = \left(y_1+x_1
z_1\right)\F_{2}\ ,\ \F_{2} z_1 = z_1 \F_{2} \ ; \F_{2} x_2
= x_2 \F_{2}\ ,\ \F_{2} y_2 = y_2 \F_{2}
$$
\end{lem}

\begin{lem}
The defining equation~(\ref{3F3}) for the operator $\F_3$
is equivalent to the system of equations
\begin{equation}
\F_{3} \left[\mathrm{L}_{1}(u_1,u_2,u_3)+
\mathrm{L}_{2}(v_1,v_2,v_3)\right]=
\left[\mathrm{L}_{1}(u_1,u_2,v_3)+
\mathrm{L}_{2}(v_1,v_2,u_3)\right] \F_{3} \label{3F3def}
\end{equation}
$$
\F_{3} x_2 = x_2 \F_{3}\ ,\ \F_{3} y_2 = y_2 \F_{3}\ ,\
\F_{3} z_2 = z_2 \F_{3}
$$
\be\F_3\cdot\left(\dd_{x_1}-z_2\dd_{y_1}\right) =
\left(\dd_{x_1}-z_2\dd_{y_1}\right)\cdot\F_3
\label{3F3last}
\end{equation}
\end{lem}
The relations in the first line are simply the rules of
commutation of $\F$-operators with $s\ell(3)$-generators
written in a compact form. In explicit notations we have
for $\Lambda_1 = (m_1, n_1)$ and $\Lambda_2 = (m_2, n_2)$
$$
\F : V_{\Lambda_1}\otimes V_{\Lambda_2} \to
V_{\Lambda^{\prime}_1}\otimes V_{\Lambda^{\prime}_2}
$$
$$
\F_1\ :\ \Lambda^{\prime}_1 = (m_1 , n_1+\xi_1)\ ;\
\Lambda^{\prime}_2 = (m_2 , n_2-\xi_1) \ ;\ \xi_1 = u_1-v_1
$$
$$
\F_2\ :\ \Lambda^{\prime}_1 = (m_1+\xi_2 , n_1-\xi_2)\ ;\
\Lambda^{\prime}_2 = (m_2-\xi_2 , n_2+\xi_2) \ ;\ \xi_2 =
u_2-v_2
$$
$$
\F_3\ :\ \Lambda^{\prime}_1 = (m_1-\xi_3 , n_1)\ ;\
\Lambda^{\prime}_2 = (m_2+\xi_3 , n_2) \ ;\ \xi_3 = u_3-v_3
$$
The $s\ell(3)$-invariance of $\R$-matrix follows directly
from the properties of $\F$-operators so that the general
R-matrix $\R^{-1}_{\Lambda_1 \Lambda_2}(v-u) =
\P_{12}\check{\R}(u;v)$ is automatically
$s\ell(3)$-invariant.

{\bf Proof.} We shall consider the operators $\F_3$ and
$\F_2$. All calculations for the operator $\F_1$ are very
similar to the $\F_3$-case. Now we are going to the proof
of equivalence of defining equation~(\ref{3F3}) to the
system~(\ref{3F3def}) and derivation of explicit formula
for the operator $\F_3$ . First we show that the
system~(\ref{3F3def}) is the direct consequence of the
eq.~(\ref{3F3}). Let us make the shift $u_k\to u_k+\lambda\
,\ v_1\to v_1+\mu\ ,\ v_2\to v_2+\nu\ ,\ v_3\to
v_3+\lambda$ in the defining equation ~(\ref{3F3}).The
$\F$-operator is invariant under this shift and
$\mathrm{L}$-operators transform as follows
$$
\mathrm{L}_1\to \mathrm{L}_1+\lambda\cdot\II \ ;\
\mathrm{L}_2\to \mathrm{L}_2 +
\lambda\cdot\II+(\mu-\lambda)\left(\begin{array}{ccc}
1 & 0 & 0\\
-x_2 & 0 & 0\\
-y_2 & 0 & 0
\end{array}\right) + (\nu-\lambda)\left(\begin{array}{ccc}
0 & 0 & 0\\
x_2 & 1 & 0\\
-x_2 z_2 & -z_2 & 0
\end{array}\right)
$$
The obtained after all equation is valid for arbitrary
parameters $\lambda$ , $\mu$ and $\nu$ and as consequence
we derive the system~(\ref{3F3def}) and
equation~(\ref{3F3last}). Next we show that from the system
of equations~(\ref{3F3def}),~(\ref{3F3last}) follows
eq.~(\ref{3F3}). This will be almost evident if we rewrite
these equations in equivalent form using the
$s\ell(3)$-invariance of the $\mathrm{L}$-operator and the
commutativity of $\F_3$ and $x_2,y_2,z_2$. We substitute
the factorized representation~(\ref{factor3}) for the
operator $\mathrm{L}_2$ in the equation~(\ref{3F3}) for the
operator $\F_{3}$ ($\mathrm{D}_{x_2}=\dd_{x_2}
-z_2\dd_{y_2}$)
$$
\F_{3} \mathrm{L}_1(u_1,u_2,u_3)\mathbf{M}
   \left(\begin{array}{ccc}
v_1 & \mathrm{D}_{x_2} & \dd_{y_2} \\
0 & v_2 & \dd_{z_2}\\
0 & 0 & v_3
   \end{array}\right)
\mathbf{M}^{-1} = \mathrm{L}_1(u_1,u_2,v_3)\mathbf{M}
   \left(\begin{array}{ccc}
v_1 & \mathrm{D}_{x_2} & \dd_{y_2} \\
0 & v_2 & \dd_{z_2}\\
0 & 0 & u_3
   \end{array}\right)
\mathbf{M}^{-1}\F_{3}
$$
and perform the similarity transformation
$\mathbf{M}^{-1}\cdots\mathbf{M}$ of this matrix equation
using the commutativity $\F_{3}$ and $x_2, y_2,z_2$. Then
using the $s\ell(3)$-invariance of
$\mathrm{L}$-operator~(\ref{sl3})
$$
\mathbf{M}^{-1}\cdot\mathrm{L}_1\cdot\mathbf{M} =
\S^{-1}\cdot \mathrm{L}_1\cdot \S\ ;\ \S =
\mathrm{e}^{z_2(\dd_{z_1}-x_1\dd_{y_1})}
\cdot\mathrm{e}^{y_2\dd_{y_1}}\cdot
\mathrm{e}^{x_2\dd_{x_1}}
\ ;\ \mathbf{M} = \left(\begin{array}{ccc}
1 & 0 & 0\\
-x_2 & 1 & 0\\
-y_2 & -z_2 & 1
   \end{array}\right)
$$
we derive the equation for the transformed operator
$\mathbf{r} = \S\cdot\F_3\cdot\S^{-1}$ \be \mathbf{r}\cdot
\mathrm{L}_1(u_1,u_2,u_3)\mathbf{L}(v_1,v_2,v_3) =
\mathrm{L}_1(u_1,u_2,v_3)\mathbf{L}(v_1,v_2,u_3)
\cdot\mathbf{r} \label{3defF} \ee where
$$
\mathbf{L}(v_1,v_2,v_3)\equiv \S
\cdot\left(\begin{array}{ccc}
v_1 & \mathrm{D}_{x_2} & \dd_{y_2} \\
0 & v_2 & \dd_{z_2}\\
0 & 0 & v_3
   \end{array}\right)\cdot\S^{-1} =
\left(\begin{array}{ccc}
v_1 & \mathrm{D}_{x_2}-\dd_{x_1}& \dd_{y_2}-\dd_{y_1}\\
0 & v_2 & \dd_{z_2}-\mathrm{D}_{z_1}\\
0 & 0 & v_3\end{array}\right)\ ;\
\mathrm{D}_{z_1}=\dd_{z_1} - x_1\dd_{y_1}
$$
To derive the system of equations which is equivalent to
the system~(\ref{3F3def}),~(\ref{3F3last}) written in terms
of $\mathbf{r}$ we repeat the same trick with the shift of
parameters and obtain the system of equations \be
\mathbf{r}\cdot
\left[\mathrm{L}_1(u_1,u_2,u_3)+\mathbf{L}(v_1,v_2,v_3)\right]
= \left[\mathrm{L}_1(u_1,u_2,v_3)+\mathbf{L}(v_1,v_2,u_3)
\right]\cdot\mathbf{r}\label{31+1} \ee
\begin{equation}
\mathbf{r}\cdot \mathrm{L}_1(u_1,u_2,u_3)
   \left(\begin{array}{ccc}
1 & 0 \\
0 & 1 \\
0 & 0
   \end{array}\right)
= \mathrm{L}_1(u_1,u_2,v_3)
   \left(\begin{array}{ccc}
1 & 0 \\
0 & 1 \\
0 & 0
   \end{array}\right)\mathbf{r}
\label{311}
\end{equation}
It is evident that all equations of the system~(\ref{311})
contained in the equation~(\ref{31+1}) except only one
$(12)$-equation $\mathbf{r} \dd_{x_1} =
\dd_{x_1}\mathbf{r}$. We use the system of equation \be
\mathbf{r}\cdot
\left[\mathrm{L}_1(u_1,u_2,u_3)+\mathbf{L}(v_1,v_2,v_3)\right]
= \left[\mathrm{L}_1(u_1,u_2,v_3)+\mathbf{L}(v_1,v_2,u_3)
\right]\cdot\mathbf{r}\ ;\  \mathbf{r} \dd_{x_1} =
\dd_{x_1}\mathbf{r} \label{3DefF} \ee as defining system
for operator $\mathbf{r}$. It is the system of equations
~(\ref{3F3def}),~(\ref{3F3last}) written in terms of
operator $\mathbf{r}$. Returning to the
system~(\ref{3defF})(which is~(\ref{3F3}) written in terms
of $\mathbf{r}$) we note that it is possible to factorize
the matrix $diag(v_1\ ; v_2\ ;\ 1)$ from the right
$$
\mathbf{r}\cdot
\mathrm{L}_1(u_1,u_2,u_3)\left(\begin{array}{ccc}
1 & \frac{\mathrm{D}_{x_2}-\dd_{x_1}}{v_2}& \dd_{y_2}-\dd_{y_1}\\
0 & 1 & \dd_{z_2}-\mathrm{D}_{z_1}\\
0 & 0 & v_3\end{array}\right) =
\mathrm{L}_1(u_1,u_2,v_3)
\left(\begin{array}{ccc}
1 & \frac{\mathrm{D}_{x_2}-\dd_{x_1}}{v_2}& \dd_{y_2}-\dd_{y_1}\\
0 & 1 & \dd_{z_2}-\mathrm{D}_{z_1}\\
0 & 0 & u_3\end{array}\right)\mathbf{r}
$$
In comparison with~(\ref{3DefF}) there are three new
equations only
$$
\mathbf{r}\cdot \mathrm{L}_1(u_1,u_2,u_3)
   \left(\begin{array}{ccc}
-\dd_{y_1}\\
-\dd_{z_1}+x_1\dd_{y_1}\\
v_3
   \end{array}\right)
= \mathrm{L}_1(u_1,u_2,v_3)
\left(\begin{array}{ccc}
-\dd_{y_1}\\
-\dd_{z_1}+x_1\dd_{y_1}\\
u_3
   \end{array}\right)
   \mathbf{r}
$$
Indeed the system~(\ref{3DefF}) contains the equations
$[\mathbf{r} , \dd_{x_2}]=[\mathbf{r} ,\partial_{y_2}]
=[\mathbf{r} ,\partial_{z_2}]=[\mathbf{r}, \dd_{x_1}]=0$
and by conditions~(\ref{311}) we obtain only three new
equations. It is possible to show that these equations
follow from the system~(\ref{3DefF}) so that it remains to
find the solution of the system of equations~(\ref{3DefF}).
First of all $[\mathbf{r} , x_2]=[\mathbf{r} , y_2]=
[\mathbf{r} , z_2]=[\mathbf{r} , \dd_{x_2}]=[\mathbf{r}
,\partial_{y_2}]=[\mathbf{r} ,\partial_{z_2}]=0$ and
therefore the operator $\mathbf{r}$ depends on the
variables $x_1,y_1,z_1$ only. There are six equations and
for simplicity we change $z_1,y_1, z_1 \to x,y,z$
\begin{equation}
\mathbf{r}\left(x\dd_x +y\dd_{y}\right) = \left(x\dd_x
+y\dd_{y}\right)\mathbf{r}\ ;\  \mathbf{r}\dd_x = \dd_x
\mathbf{r}\ ;\ \mathbf{r}\left(-x\dd_x +z\dd_{z}\right) =
\left(-x\dd_x +z\dd_{z}\right)\mathbf{r} \label{331}
\end{equation}
\begin{equation}
\mathbf{r}\left(y\dd_x +z^2\dd_{z}+(u_2-u_3+1)z\right) =
\left(y\dd_x +z^2\dd_{z}+(u_2-v_3+1)z\right)\mathbf{r}
\label{332}
\end{equation}
\begin{equation}
\mathbf{r}\left(x^2\dd_x
+xy\dd_{y}-xz\dd_z-y\dd_z+(u_1-u_2+1)x\right) = \label{333}
\end{equation}
$$
=\left(x^2\dd_x+xy\dd_{y}-xz\dd_z-y\dd_z+
(u_1-u_2+1)x\right)\mathbf{r}
$$
\begin{equation}
\mathbf{r}\left(x\left(y\dd_x+z^2\dd_z+(u_2-u_3+1)z\right)
+y\left(y\dd_y+z\dd_z+u_1-u_3+2\right)\right) = \label{334}
\end{equation}
$$
=\left(x\left(y\dd_x+z^2\dd_z+(u_2-v_3+1)z\right)
+y\left(y\dd_y+z\dd_z+u_1-v_3+2\right)\right)\mathbf{r}
$$
We look the general solution in the form
$$
\mathbf{r} = \mathbf{a}[z\dd_z]\cdot
\mathrm{e}^{\frac{y}{z}\dd_x}\cdot
\mathbf{b}[y\dd_y]\cdot\mathrm{e}^{-\frac{y}{z}\dd_x}
\cdot\mathbf{c}[z\dd_z]
$$
where
$$
\mathbf{c}[z\dd_z]=
\frac{\Gamma(z\dd_z+1)}{\Gamma(z\dd_z+u_2-u_3+1)}.
$$
It is the evident solution of the equations~(\ref{331}).
The equations~(\ref{332})-(\ref{334}) lead to the
recurrence relations for the functions $\mathbf{a}[z\dd_z]$
and $\mathbf{b}[y\dd_y]$
$$
\mathbf{a}[z\dd_z+1]\cdot(z\dd_z+1) =
(z\dd_z+u_2-v_3+1)\cdot\mathbf{a}[z\dd_z]
$$
$$
\mathbf{b}[y\dd_y+1]\cdot(y\dd_y+u_1-u_3+1) =
(y\dd_y+u_1-v_3+1)\cdot\mathbf{b}[y\dd_y]
$$
The solution of these equations is
$$
\mathbf{a}[z\dd_z]
=\frac{\Gamma(z\dd_z+u_2-v_3+1)}{\Gamma(z\dd_z+1)}\ ;\
\mathbf{b}[y\dd_y] =
\frac{\Gamma(y\dd_y+u_1-v_3+1)}{\Gamma(y\dd_y+u_1-u_3+1)}.
$$
Collecting everything together we obtain the expression for
the operator $\F_3$ from the Proposition. All calculations
for the operator $\F_1$ are very similar and finally one
arrives to the system which coincides with the
system~(\ref{331}),~(\ref{334}) after the change of
variables and parameters. In a such way one obtains the
expression for the operator $\F_1$ from the Proposition. It
remains to prove the equivalence of defining
equation~(\ref{3F2}) to the system~(\ref{3F2def}) and
derive the explicit formula for the operator $\F_2$. First
we show that the system~(\ref{3F2def}) is the direct
consequence of the eq.~(\ref{3F2}). Let us make the shift
$u_1\to u_1+\lambda\ ,\ u_2\to u_2+\lambda\ ,\ u_3\to
u_3+\mu\ ,\ v_1\to v_1+\nu\ ,\ v_2\to v_2+\lambda\ ,\
 v_3\to v_3+\lambda$ in the defining equation ~(\ref{3F2})
for the operator $\F_2$.The $\F$-operator is invariant
under this shift and $\mathrm{L}$-operators transform as
follows
$$
\mathrm{L}_1\to \mathrm{L}_1+ \lambda\cdot\II
+(\mu-\lambda)\left(\begin{array}{ccc}
0 & 0 & 0\\
0 & 0 & 0\\
y_1+x_1 z_1 & z_1 & 1
\end{array}\right)\ ;\
\mathrm{L}_2 \to \mathrm{L}_2 + \lambda\cdot\II+
(\nu-\lambda)\left(\begin{array}{ccc}
1 & 0 & 0\\
-x_2 & \lambda & 0\\
y_2 & 0 &0\end{array}\right)
$$
The obtained after all equation has to be satisfied for
arbitrary $\lambda$ , $\mu$ and $\nu$ and as consequence we
derive the system~(\ref{3F2def}). Next we show that from
the systems of equations~(\ref{3F2def}) follows
eq.~(\ref{3F2}). This will be almost evident if we rewrite
these equations in equivalent form using the
$s\ell(3)$-invariance of the $\mathrm{L}$-operator and the
commutativity of $\F_2$ and $z_1 , y_1+x_1z_1 , x_2 , y_2$.
We start from the equation~(\ref{3F2}) and make the two
similarity transformations of the defining
equation~(\ref{3F2}) using simple matrices which commute
with operator $\F_2$. After all these transformations the
defining equation~(\ref{3F2}) for the $\F_2$-operator in
factorized form looks as follows
$$
\F_2\cdot\mathbf{l}_1(u_1,u_2,u_3)\cdot
\mathbf{M}\cdot\mathbf{l}_2(v_1,v_2,v_3) =
\mathbf{l}_1(u_1,v_2,u_3)\cdot
\mathbf{M}\cdot\mathbf{l}_2(v_1,u_2,v_3) \cdot\F_2\ ;\
\mathbf{M}\equiv \left(\begin{array}{ccc}
1 & 0 & 0 \\
-x_2 & 1 & 0 \\
y_{12} + z_1 x_{12} & z_1 & 1
   \end{array}\right)
$$
$$
\mathbf{l}_1(u_1,u_2,u_3)\equiv\left(\begin{array}{ccc}
1 & 0 & 0 \\
-x_1 & 1 & 0 \\
0 & 0 & 1
   \end{array}\right)\left(\begin{array}{ccc}
u_1 & \dd_{x_1}-z_1\dd_{y_1} & \dd_{y_1} \\
0 & u_2 & \dd_{z_1} \\
0 & 0 & u_3
   \end{array}\right)\left(\begin{array}{ccc}
1 & 0 & 0 \\
x_1 & 1 & 0 \\
0 & 0 & 1
   \end{array}\right)
$$
$$
\mathbf{l}_2(v_1,v_2,v_3)\equiv \left(\begin{array}{ccc}
1 & 0 & 0 \\
0 & 1 & 0 \\
0 & -z_2 & 1
   \end{array}\right)\left(\begin{array}{ccc}
v_1 & \dd_{x_2}-z_2\dd_{y_2}  & \dd_{y_2} \\
0 & v_2 & \dd_{z_2} \\
0 & 0 & v_3
   \end{array}\right)\left(\begin{array}{ccc}
1 & 0 & 0 \\
0 & 1 & 0 \\
0 & z_2 & 1
   \end{array}\right)
$$
Next step we rewrite this equation in terms of the
transformed operator $\mathbf{r}$
$$
\F_2 = \S^{-1} \cdot\mathbf{r}\cdot\S\ ;\ \S =
\mathrm{e}^{\left(y_2-x_1z_1\right)\partial_{y_1}}\cdot
\mathrm{e}^{z_1\partial_{z_2}}\cdot
\mathrm{e}^{x_2\partial_{x_1}}
$$
\be \mathbf{r}\cdot\mathbf{L}_1(u_1,u_2,u_3)\mathbf{m}
\mathbf{L}_2(v_1,v_2,v_3) =
\mathbf{L}_1(u_1,v_2,u_3)\mathbf{m}
\mathbf{L}_2(v_1,u_2,v_3) \cdot\mathbf{r}\ ;\ \mathbf{m}=
\left(\begin{array}{ccc}
1 & 0 & 0 \\
0 & 1 & 0 \\
y_{1} & 0 & 1
   \end{array}\right)\label{3defF2}
\ee where
$$
\mathbf{L}_1(u_1,u_2,u_3) = \left(\begin{array}{ccc}
u_1+1+x_1\dd_{x_1} & \dd_{x_1} & \dd_{y_1} \\
-x_1^2\dd_{x_1}+(u_2-u_1-1)x_1 & u_2-x_1\dd_{x_1} &
\dd_{z_1}-\dd_{z_2} \\
0 & 0 & u_3
   \end{array}\right)
$$
$$
\mathbf{L}_2(v_1,v_2,v_3) = \left(\begin{array}{ccc}
v_1 & \dd_{x_2}-\dd_{x_1}-z_1\dd_{y_2} & \dd_{y_2}-\dd_{y_1} \\
0 & v_2+1+z_2\dd_{z_2} & \dd_{z_2} \\
0 & -z_2^2\dd_{z_2} +(v_3-v_2-1)z_2 & v_3-z_2\dd_{z_2}
   \end{array}\right)
$$
To derive the system of equations which is equivalent to
the system~(\ref{3F2def}) written in terms of $\mathbf{r}$
we repeat the same trick with the shift of parameters and
obtain \be
\mathbf{r}\cdot\left[\mathbf{L}_1(u_1,u_2,u_3)\mathbf{m}+
\mathbf{m}\mathbf{L}_2(v_1,v_2,v_3)\right] =
\left[\mathbf{L}_1(u_1,v_2,u_3)\mathbf{m}+
\mathbf{m}\mathbf{L}_2(v_1,u_2,v_3)\right]\cdot\mathbf{r}
\label{32+2} \ee This system results in a simple equations.
First of all we have $[\mathbf{r}, y_1]=[\mathbf{r},
z_1]=[\mathbf{r}, x_2]=[\mathbf{r}, y_2]=0$ from the very
beginning and the system contain the equations \be
[\mathbf{r}, \dd_{y_2}] = [\mathbf{r}, \dd_{z_1}] =
[\mathbf{r}, \dd_{x_2}-z_1\dd_{y_2}] = 0 \label{321}\ee \be
\mathbf{r}\left(x_1\dd_{x_1}+y_1\dd_{y_1}\right) =
\left(x_1\dd_{x_1}+y_1\dd_{y_1}\right)\mathbf{r}\ ;\
\mathbf{r}\left(x_1\dd_{x_1}-z_2\dd_{z_2}\right) =
\left(x_1\dd_{x_1}-z_2\dd_{z_2}\right)\mathbf{r}
\label{322} \ee \be \mathbf{r}\left(-x_1^2\dd_{x_1} +
(u_2-u_1-1)x_1 -y_1\dd_{z_2}\right) = \left(-x_1^2\dd_{x_1}
+ (v_2-u_1-1)x_1 -y_1\dd_{z_2}\right)\mathbf{r} \label{323}
\ee \be \mathbf{r}\left(-z_2^2\dd_{z_2} + (v_3-v_2-1)z_2
-y_1\dd_{x_1}\right) = \left(-z_2^2\dd_{z_2} +
(v_3-u_2-1)z_2 -y_1\dd_{x_1}\right)\mathbf{r} \label{324}
\ee Returning to the system~(\ref{3defF2}) which is
equivalent to the system~(\ref{3F2}) written in terms of
$\mathbf{r}$ we note that it is possible to factorize the
matrix $diag(v_1 ; 1 ; 1)$ from the right and the matrix
$diag(1 ; 1 ; u_3)$ from the left so that the parameters
$v_1$ and $u_3$ disappear from equation. The equivalence
obtained system of equations to the
system~(\ref{321})-(\ref{324}) can be proven by
straightforward calculations. Now we are going to the
solution of the defining system of equations. We look the
general solution in the form
$$
\mathbf{r} = \mathbf{a}[z_2\dd_{z_2}]\cdot
\mathrm{e}^{\frac{y_1}{z_2}\dd_{x_1}}\cdot
\mathbf{b}[x_1\dd_{x_1}]\cdot\mathrm{e}^{-\frac{y_1}{z_2}\dd_{x_1}}
\cdot\mathbf{c}[z_2\dd_{z_2}]
$$
where
$$
\mathbf{c}[z_2\dd_{z_2}]
=\frac{\Gamma(z_2\dd_{z_2}+1)}{\Gamma(z_2\dd_{z_2}+v_2-v_3+1)}.
$$
It is the evident solution of the
equations~(\ref{321}),~(\ref{322}). The
equations~(\ref{333}) and~(\ref{334}) lead to the
recurrence relations for the functions
$\mathbf{a}[z_2\dd_{z_2}]$ and $\mathbf{b}[x_1\dd_{x_1}]$
$$
\mathbf{a}[z_2\dd_{z_2}+1]\cdot(z_2\dd_{z_2}+1) =
(z_2\dd_{z_2}+u_2-v_3+1)\cdot\mathbf{a}[z_2\dd_{z_2}]
$$
$$
\mathbf{b}[x_1\dd_{x_1}+1]\cdot(x_1\dd_{x_1}+u_1-u_2+1) =
(x_1\dd_{x_1}+u_1-v_2+1)\cdot\mathbf{b}[x_1\dd_{x_1}]
$$
The solution of these equations is
$$
\mathbf{a}[z_2\dd_{z_2}]
=\frac{\Gamma(z_2\dd_{z_2}+u_2-v_3+1)}{\Gamma(z_2\dd_{z_2}+1)}\
;\ \mathbf{b}[x_1\dd_{x_1}] =
\frac{\Gamma(x_1\dd_{x_1}+u_1-v_2+1)}{\Gamma(x_1\dd_{x_1}+u_1-u_2+1)}
$$
and collecting everything together we obtain the expression
for the operator $\F_2$ from the proposition.

\section{Conclusions}

We have shown that the general R-matrix can be represented
as the product of the simple "building blocks" --
$\F$-operators. In the present paper we have demonstrated
how this factorization arises in the simplest situations of
the symmetry algebra $s\ell(2)$ and $s\ell(3)$. As a
byproduct we derived useful representation for the
$s\ell(2)$-invariant R-matrix
$$ \R_{\ell_1\ell_2}(u) =
\P_{12}\RR_{1}(u_1|v_1,u_2)\RR_{2}(u_1,u_2|v_{2})\sim
\P_{12}\cdot \frac{\Gamma(z_{21}\dd_2+2\ell_1)}
{\Gamma(z_{21}\dd_2+\ell_1+\ell_2-u)}
\cdot\frac{\Gamma(z_{12}\dd_1+\ell_1+\ell_2+u)}
{\Gamma(z_{12}\dd_1+2\ell_1)}
$$
It is possible to represent the general
$s\ell(2)$-invariant R-matrix in many equivalent forms. The
spectral decomposition of R-matrix~(\ref{spec}) was
obtained in the paper~\cite{KRS}. In functional
representation R-matrix is some integral operator acting on
the space of polynomials. The symbol and kernel of this
integral operator was calculated in~\cite{Skl} and
in~\cite{DKK} correspondingly. Here we obtain the
representation for the R-matrix in the operator form. The
similar operator expression for the Hamiltonian of the XXX
spin chain was used in~\cite{L,FK,KK}.

In the $s\ell(3)$ case we have derived the explicit
expression for the general R-matrix in the factorized form
$$
\R^{-1}_{\Lambda_1\Lambda_2}(v-u) \sim \P_{12}
\F_1(u_1;v_1,u_2,u_3)
\F_2(u_1,u_2;v_2,u_3)\F_3(u_1,u_2,u_3;v_3).
$$
The $\F$-operators in many respects are very similar to the
general R-matrix but are much more simpler. In some sense
the R-matrix is composite object and $\F$-operators are
elementary building blocks. Using the simple pictures it is
possible to derive the general system of defining relations
for the $\F$-operators. First, there exist nontrivial
commutation relations for the $\F$-operators acting in
different spaces. Second, there are three term relations
for the $\F$-operators which play the same role as
Yang-Baxter relation for the R-matrix. Using defining
relations for $\RR$-operators it is possible to derive all
relations for the general R-matrix including the
Yang-Baxter equation. All this will be discussed elsewhere.

It seems that all these results can be generalized to the
symmetry algebra $s\ell(n)$ where should exists the
factorization of the general R-matrix on the product of n
simple $\F$-operators. In the second part of this work we
shall show that the same factorization take place for the
general rational solution of the Yang-Baxter equation with
the supersymmetry algebra $s\ell(2|1)$.

\section{Acknowledgments}

I would like to thank R.Kirschner, G.Korchemsky, P.Kulish,
A.Manashov, E.Sklyanin and~V.Tarasov for the stimulating
discussions and critical remarks on the different stages of
this work. This work was supported by the grant 03-01-00837
of the Russian Foundation for Fundamental Research.

\end{document}